\documentclass[12pt]{amsart}
\usepackage{graphicx}
\usepackage{amssymb}
\usepackage{epstopdf}
\usepackage{comment}
\usepackage{breqn}
\usepackage{color}

\newtheorem{theorem}{Theorem}[section]
\newtheorem{prop}[theorem]{Proposition}
\newtheorem{lemma}[theorem]{Lemma}
\newtheorem{cor}[theorem]{Corollary}
\newtheorem{definition}[theorem]{Definition}

\newcommand{\btheorem}{\begin{theorem}}
\newcommand{\etheorem}{\end{theorem}}
\newcommand{\bprop}{\begin{prop}}
\newcommand{\eprop}{\end{prop}}
\newcommand{\blemma}{\begin{lemma}}
\newcommand{\elemma}{\end{lemma}}
\newcommand{\bcor}{\begin{cor}}
\newcommand{\ecor}{\end{cor}}

\theoremstyle{remark}
\newtheorem{remark}[theorem]{Remark}

\newcommand{\bg}{{\mathfrak g }}

\newcommand{\bl}{{\mathfrak l }}

\newcommand{\bq}{{\mathfrak q }}

\newcommand{\bh}{{\mathfrak h }}
\newcommand{\bt}{{\mathfrak t }}

\newcommand{\bR}{{\mathbb R}}
\newcommand{\bN}{{\mathbb N}}
\newcommand{\bC}{{\mathbb C}}
\newcommand{\bZ}{{\mathbb Z}}

\newcommand{\bL}{{\mathbb L}}

\newcommand{\C}{{\mathbb C}}
\newcommand{\R}{{\mathbb R}}
\newcommand{\Z}{{\mathbb Z}}
\newcommand{\A}{{\mathbf A}}
\newcommand{\N}{{\mathbb N}}
\newcommand{\IH}{{\mathbb H}}

\usepackage[parfill]{parskip}

\title{Some branching laws \linebreak
for symmetric spaces.}
\author{Bent \O{}rsted \\ Birgit Speh}

\address{
B. \O{}rsted\\               
Department of Mathematics \\
Aarhus University \\ 
Ny Munkegade 118 \\
Building 1530, 431 \\
8000 Aarhus C \\
Denmark \\
orsted@math.au.dk                
}

\address{
B. Speh \\               
Department of Mathematics \\ Malott Hall, Cornell University \\ Ithaca, NY 14853 \\ 
United States\\           
bes12@cornell.edu                
}

\begin{document}


\maketitle

\begin{abstract}In this paper we consider the unitary symmetric spaces of the form
$X = U(p,q)/U(1)U(p,q-1)$ and their discrete series representations.
Inspired by the work of A.~Venkatesh and Y.~Sakellaridis on $L$-groups  of $p$-adic spherical   spaces \cite{SV} we formulate and prove natural relative branching laws for the restriction of such discrete series representations to smaller groups of the same type and corresponding unitary symmetric spaces; here all the spaces involved have the same $L$-group. We think of this as steps towards 
understanding Gan Gross Prasad
conjectures for unitary symmetric spaces. Indeed, using  period integrals, some results on the Laplacian of spheres, \cite{STV} and results by T.~Kobayashi \cite{K} we prove  an  analogue of these conjectures. 


\end{abstract}

\

\begin{center}{\em

It is pleasure to dedicated this article to Toshiyuki Kobayashi whose work is fundamental to the understanding of the restriction of representations of reductive groups to subgroups}

\end{center}


\section{introduction}
In this paper we shall study in detail some branching laws for unitary representations
corresponding to a pair of Lie groups $G' \subset G$. We consider a unitary irreducible
representation $\Pi$ of $G$ in a Hilbert space $\mathcal H$ and its restriction to the subgroup $G'$. The
branching law gives the decomposition of $\mathcal H$ as a direct integral of unitary
irreducible representations $\pi$ of $G'$, possibly with multiplicities. This is the analogue
of the spectral decomposition of a self-adjoint linear operator on $\mathcal H$, and of
special importance is the discrete spectrum. In broad terms we shall consider families of
representations $\Pi$ as irreducible invariant subspaces of $L^2(G/H)$,
the so-called discrete series of a symmetric space $G/H$ (where $H$ is non-compact).  Their Langlands parameters are determined by H. Schlichtkrull, \cite{Sch2}. They can also be parametrized using the $L$-group of the symmetric space; see \cite{MR} and 
\cite{SV}.
The motivation for this work is to obtain branching laws for discrete series representations of symmetric spaces, which are similar to the branching laws of the
Gross-Prasad conjectures for discrete series representations of the unitary group, see \cite{GP}. We are in particular interested in the representations $\pi$ in the branching
law that belong to a similar family  coming from the discrete
spectrum of a symmetric subspace $G'/H' \subset G/H$. 

There are essentially three issues 
to clarify in this context, namely
\begin{itemize}
\item The  representations $\Pi$  appear as irreducible subspaces of
the left regular representation $L^2(G/H)$, concretely for $G = U(p,q)$ and $H$ is either $ U(1) \times U(p-1,q)$
or $ U(1) \times U(p,q-1)$.
They  are obtained by  cohomological
induction from a character of a  $\theta$-stable parabolic subgroup. Considered  as members of a Arthur-Vogan packet they also have an associated epsilon character on $\bZ_2^{2}$ \cite{ABV}.
 \item The subgroup $G'$ in the branching law that we study is $G' = U(p-1,q)$; here
the choice of the imbedding is important, and the parameters for the representations
$\Pi$ are again the cohomological parameter, i.e a character of a $\theta$-stable parabolic subgroup and  an epsilon character using  the data
corresponding to a discrete series representation for a subspace $G'/H'$.
\item There are essentially two qualitatively different situations that we have to consider to determine the
branching corresponding to $G' \subset G$. In one case  there is both a continuous spectrum
(where purely algebraic methods are difficult to apply) and where we focus  on the discrete spectrum. In the other case the spectrum is purely
discrete (and even admissible, i.e. the decomposition is a direct sum of irreducible representations with
finite multiplicity). 
\end{itemize}

Note that the representations $\Pi $ are typically not discrete series \linebreak representations of $G$, hence they are not covered by the prediction of 
 the Gross-Prasad conjectures. At the same time these representations occur in natural families as a family of representations  in the discrete spectrum of  spherical spaces $G/H$.

We shall state the parameters of the representations in a way that makes the branching law
look simple and natural, and use in the proof the idea of branching in stages corresponding to a sequence
$G' \subset G \subset G_1$ for a convenient group $G_1$; we shall also use certain
period integrals to illustrate the explicit coupling between a $\Pi$ and a $\pi$ in the
branching law (representations of $G$ resp. $G'$). 

To explain in a little more detail our ideas, consider first compact groups.

  A classical problem in the representation theory of a connected compact group Lie group $G$ is the restriction of a representation of $G$ to a (usually connected) subgroup $G'$. The restriction is a direct sum of irreducible representations of the subgroup.  In some cases it is easy (though it might in practice { involve }
complicated combinatorics) to characterize the representations  of $G'$ which occur in the restriction with positive multiplicities using  branching rules - which are often stated and studied in introductory textbooks. For example the classical restriction  of irreducible representations of the unitary group $U(n)$  to $U(n-1)$ can be found 
also on Wikipedia, see  
https://en.wikipedia.org/wiki/Restricted-representation. For more general pairs of groups the branching rules were first formulated by B. Kostant in an unpublished note generalizing the methods in  \cite{Ko}. The multiplicities are very complicated combinatorial expressions \cite{V}.

 In {\it relative branching} we consider the special case of families of representations with a given geometric realization for a compact group $G$ and ask for nonzero multiplicities of
representations of the connected subgroup  $G'$ in the restriction which have a similar geometric realization. For example we might
have a symmetric subgroup $H \subset G$ and the corresponding family of spherical representations $V$ of $G$, which 
 are realized in the  functions on $X = G/H$. Now we consider the similar family for $G'$. We consider 
$H' = G' \cap H$ (assume $G'$ is stable under the symmetry defining $H$) and $X' = G'/H'$, and 
representations of $G'$ are realized in the functions on $X'$. We want to find the multiplicities
\[  \mbox{dim Hom}_{G'}( V \otimes (V') ^*,\bC)\]  of the $H'$-spherical
representations $V'$ of $G'$ in a given $H$-spherical representation $V$ of $G$. A geometric  natural way
to do this is to study the {\it period integral}
$$I = \int_{X'} \psi(x') \psi'(x')  dx'$$
i.e. the natural pairing of the functions $\psi \in V \subset C^\infty(G/H) $ and $\psi' \in (V')^*\subset C^\infty(G'/H') $.
 In particular we may ask for
the convergence, and the non-vanishing (and exact value) of $I$ for  $H$-invariant spherical functions $\psi$, respectively $H'$-invariant spherical functions $\psi'.$
  Even  in the compact case
 it is an interesting
conjecture that the period  integral  for  spherical functions  of finite dimensional representations is nonzero if and only if the $H'$-spherical $V'$ occurs in the
branching from $V$. These questions also make sense for non-compact groups,
as we shall see and study 
in this paper.



 
 
For unitary representations of
 noncompact reductive Lie groups $G$ the situation is more complicated since the restriction of an irreducible representation $\Pi$ to a subgroup $G'$ usually does not decompose into a direct sum of irreducible unitary representations and we may have both a continuous and also a discrete spectrum in the direct
integral decomposition (over the unitary dual of $G'$) of the restriction. Never the less considering $G \times G/ G$ as a symmetric space and the discrete series representations as its discrete spectrum, period integrals are  important for the restriction of discrete series  representations  to large subgroups \cite{OV}.

 In series of papers starting 1991  B. Gross and D. Prasad \cite{GP} introduced conjectures for the restriction of discrete series representations of orthogonal groups to orthogonal subgroups, These  were generalized by Gan, Gross, and Prasad to the restriction of discrete series representations of unitary groups $U(p,q)$ to unitary subgroups \cite{GGP} and  were proved by  Hongyu He   for real unitary groups \cite{H}
 and by R. Beuzart-Plessis  in the $p$-adic case. \cite{BP}
 An important role in  these conjectures   is played by Vogan packets and interlacing relations between the infinitesimal characters of the discrete series representations.

The aim of this article is to
generalize this circle of ideas  as well as  the considerations for orthogonal groups in the article "Hidden Symmetries" by T. Kobayashi and one of the authors. \cite{KS}.  In this we use ideas of  A. Venkatesh and Y. Sakellaridis who define for a  spherical space $X=G/H$ a complex Levi subalgebra ${\bl}_X$ and an $L$-group \cite{SV}. In this article we
consider  several symmetric spaces $G/H$ for a unitary group $ G=U(p,q)$ and $H=G^\sigma$ for an involution $\sigma$.  We assume that  these symmetric spaces have the complex Levi subalgebra 
$\bL_X=gl(n-2,\C) \oplus \C^2$ and that the Lie algebra of its $L$-group is  $sp(2,\bC) $ (as defined  by A. Venkatesh  and Y. Sakellaridis ). The unitary symmetric spaces satisfying these assumptions are
 \[X^+= G/H^+ =U(p,q)/U(1,0)U(p-1,q) \] and \[X^-= G/H^-=U(p,q)/U(0,1)U(p,q-1).\] 
 Here  $H^+$ is  the stabilizer of the first basis vector $e_1$ and $H^- $  the stabilizer of $e_{p+1}.$
 For $X^+$ the Levi subgroup is \[L^+=U(1)U(1) U(p-2,q)\] and for $X^-$ it 
 is \[L^-=U(1)U(1)U(p,q-2).\] 
 
 We consider a pair of representations $\Pi^+$  and  $\Pi^-$  in the discrete spectrum  of $L^2(X^+)$, respectively $L^2(X^-)$, with the same infinitesimal character. These representations are cohomologically induced from a character of a $\theta $-stable parabolic with Levi subgroup
$L^+=U(1)U(1)U(p-2,q),$ respectively $L^-= U(1)U(1)U(p,q-2).$
The work of H. Schlichtkrull and E. van den Ban implies that they have multiplicity one in the discrete spectrum  of exactly one of the symmetric spaces \cite{vB}. Following Moeglin/Renard  \cite{MR} we consider an Arthur packet $A(\phi)$ which contains both of these representations and define epsilon characters. 

Let $G' =U(p-1,q)$  be fixed point group of  the basis vector $e_{p}$. We have symmetric spaces
\begin{eqnarray*} 
Y^+ & = & G'/H^+ \cap  U(p-1,q)\\
  \ & = &U(p-1,q)/U(1)U(p-2,q) 
  \end{eqnarray*} 
  and 
   \begin{eqnarray*} Y ^- & = & G'/H^- \cap U(p-1,q) \\
  \ & = &U(p-1,q)/U(1)U(p-1,q-1). 
  \end{eqnarray*}
   The symmetric spaces $Y^+$ and $Y^-$ have  the same L-group and the same complex Levi subalgebra. 
   


 We consider in this article the relative restriction of  representations $\Pi^+, \Pi^- $ in the discrete spectrum of the symmetric spaces $L^2(X^+)$ and $L^2(X^-)$  to the discrete spectrum  of the symmetric spaces
$Y^+\subset X^+$ respectively $Y^-\subset X^-$. 

\medskip
To avoid considering special cases we assume  in this article           \begin{center}   $p \ge 3$ and $q \ge 2$.  \end{center}
\medskip

 Under this assumption all the symmetric spaces  have a nonempty discrete $ L^2$--spectrum and the representations in the discrete $L^2$--spectrum are parametrized by a character of the Levi of a $\theta$ stable parabolic and they have an epsilon character.  


\medskip
The goal of this article is to find and discuss the relative branching laws for the restriction of representations on the discrete spectrum of $X^+$ and $X^-$ to the discrete spectrum of $Y^+$ and $Y^1$.  To understand the branching  we study both the analytic models of the unitary representations 
appearing in the discrete spectrum of the symmetric spaces and the corresponding periods
integrals.

Our first conclusion will be as follows:

{ Following here the conventions of the article by C.Moeglin and D. Renard \cite{MR} for the parametrization of cohomological induced representations (see \ref{sec:discrete spec} for more details). Let $a\geq \frac{p+q}{2}$ be a nonnegative integer if $p+q$ is even and a positive $1/2$ integer  if $p+q$ is odd. 
Then $\lambda^+=(a,-a, 0,\dots ,0) $,  is a parameter  for $\Pi^+$, and $\lambda^-=(0,\dots,0,a,-a)$ a parameter for $\Pi^-$ of $U(p,q)$ in the discrete spectrum of $X^+$, respectively $X^-$. Let $b$ is a non-negative integer if $p+q-1$ is even and an 1/2 integer if $p+q -1$ is odd and $b \geq \frac{p+q-1}{2}$.
 Then $(\lambda')^+  =(b,-b,\dots ,0)$,$(\lambda')^- =( 0,\dots, 0,b,-b)$  is a parameter of a representation $\pi^+$
 respectively $\pi^-$ of $G'=U(p,q-1)$ in the discrete spectrum of $Y^+$ respectively $Y^-.$ We denote the representations by $\Pi^+_a$ or $\Pi^-_a$ respectively $\pi^+_b$ or $\pi^-_b $ . } 
 
 \begin{theorem} \label{theorem:1}
  Assume that G and G'  are 
as discussed above and that the representations in the discrete spectrum of the symmetric spaces are parametrized as above. Then
  \begin{enumerate} 
  \item 
   \[ \mbox{dim Hom}_{G'} (\Pi^+ _a , \pi_b  ^+) =1 \]
  
 if and only if $a > b$ (and zero otherwise) and
 \[
    \mbox{dim Hom}_{G'} (\Pi^+ _a , \pi_b  ^-) =0 \]
 \item  
     \[ \mbox{dim Hom}_{G'} (\Pi^- _a , \pi_b  ^-) =1 \]
 if and only if $a < b$ (and zero otherwise) and
 \[ \mbox{dim Hom}_{G'} (\Pi^- _a , \pi_b  ^+) =0.\]
 
\end{enumerate}
 \end{theorem}
Note the distinction here with finitely many representations in this relative branching
in the first case and infinitely many in the second. This is (as we shall see in the proof)
the distinction between non-admissible and admissible branching.
 
 
 
  
  
In section V.3 we reformulate this theorem using interlacing patterns of infinitesimal characters and epsilon-characters.



\medskip
 Let $\lambda$ be an integral nonsingular infinitesimal character for $U(p,q)$ and $\lambda'$ an infinitesimal character of $G'=U(p,q-1)$. We denote by $\Pi^+(\lambda)$, $\Pi^-(\lambda)$ the representations in the discrete spectrum of $X^+$, respectively $X^-$  with infinitesimal character $\lambda$. 
Define 
\[A_S (\lambda) = \Pi^+(\lambda) \oplus \Pi^-(\lambda )\]

Similarly we define  for an integral regular infinitesimal character  $\lambda'$  
 for $G'= U(p,q-1)$  the representations $\pi^+(\lambda') $ and $\pi^-(\lambda')$ and define \[B_S'(\lambda')= \pi^+(\lambda') \oplus \pi^-(\lambda')\]
  Following the ideas  in \cite{GP}  we consider
\[\mbox{Hom}_{G'}(A_S(\lambda)_{G'}, B_S(\lambda')  ).\]

 The  formulation of conjectures by B, Gross and D. Prasad for discrete series representations, our relative branching laws and the observation  that only one of the integers $2a, 2b$ is even inspired the following reformulation of Theorem \ref{theorem:1} - note that the infinitesimal characters $\lambda$ and
$\lambda'$ correspond to the parameters $a$ and $b$ above.

\begin{theorem}\label{theorem:2}
Assume that  $p,q > 3$ and $p \not = q$ and that a,b satisfy the previous assumptions. Using the above notation we have



\[\mbox{dim Hom}_{G'}(A_S(\lambda)_{|G"}, B_S(\lambda')  )= 1.\]
\end{theorem}

\medskip

The paper is organized as follows. In section II. we introduce our notation. In Section III. we discuss the symmetric spaces, their $L$-groups and the representations in the discrete spectrum. In section IV the concept of relative branching is introduced and illustrated by examples. The relative branching laws are formulated in section V. as well some interpretations and a similarity to  a branching law  to   Gross-Prasad branching laws  using the epsilon characters of the representations. This is one of our key motivations,
giving the relative branching law in terms of parameters coming from the theory of spherical spaces.
 In section VI. we return to the proofs of the branching laws using some classical techniques;
  we identify $Y^+= U(p,q)/U(1)U(p-1,q-1) $  with an orbit of $G'$ on $X^+$.   Flensted-Jensen in his seminal work \cite{FJ}
constructed discrete series for some symmetric spaces $G/H$ and we shall use his realization of the representations  in our
special case $L^2(X^+)$.
We consider a geometric restriction of the representations in the discrete spectrum of $L^2(X^+) $ to  $L^2(Y^+)$  by using the restriction of the $K$-finite coefficients of the Flensted Jensen representation to  $Y^+$  and compute the corresponding  period integral.
  Some of the technical points  use recent work of T. Kobayshi \cite{K} on $O(p,q)$ as well
as work of Schlichtkrull, Trapa, and Vogan \cite{STV} on $U(p,q)$ - these are the groups
$G$ in question. In partiular in Section VI. we prove the exhaustion of the
branching law predicted by the period integrals. Note here the sequences
\begin{itemize}
\item $U(p-1,q) \subset O(2p-2, 2q) \times O(2) \subset O(2p, 2q)$
\item $U(p-1,q) \subset U(p,q) \subset O(2p,2q)$
\end{itemize}
that we shall study in Section VI. In particular one may do branching in stages along either
sequence of groups and obtain information about the first branching in the second line. 
   In  section VII.  we discuss the branching laws for $\Pi^-$ \cite{S}. The restriction is an infinite sum of irreducible representations and we determine the irreducible summands which are in the discrete spectrum of $L^2(Y^-)$.  We determine again  $K$-finite functions in $\Pi^- \subset L^2(X^-)$ which generate the  irreducible $(\bg',K')$-modules of these summands and then show that the restriction of these functions to $Y^- \subset X^-$ is non zero.   
   The  restriction of $\Pi^-$ to $G'$ was determined 
 using different techniques  and in greater generality by Yoshiki Oshima around 2014 and his results show that there are representations  in $\Pi^-_{|G'}$ which are not in the discrete spectrum of $L^2(Y^-)$.
In an appendix we discuss relative branching and  period integrals for compact Riemannian symmetric
spaces of rank one,
extending the discussion in particular to the quaternionic case.

              
    The second author  was introduced to the ideas of Y.~Sakellaridis and A.~Venkatesh about discrete series representations of spherical spaces  in a workshop organized at AIM by A.~Venkatesh. She  would also like to  thank  D.~Vogan and J. ~Adams for helpful discussions about Arthur packets, S.~Sahi for many  conversations about orthogonal polynomials and Yoshiki Oshima for sharing his manuscript about branching for a family of cohomologically induced representations of unitary groups. { Much of the joint work was done using the facilities of the virtual AIM research group "Representation theory and non commutative Operator algebras". }

.

{\bf Notation:}: $\bN = \{0,1,2,\dots ,\} $ and  $\bN_+ = \{1,2, \dots ,\}$.  


\bigskip

\section{Notation and generalities}
 
\subsection{} {
Consider the hermitian quadratic form $Q(Z,Z) $ on ${\mathbb{C}}^{p+q}$  with signature (p,q).  We choose the basis  $e_1, \dots ,e_p, e_{p+1},\dots e_{p+q}$ so that 
 the form is positive on $e_1, \dots ,e_p$ and negative on $e_p, e_{p+1},\dots e_{p+q}$ and 
\begin{equation}
\label{eqn:quad}
  Q(Z,Z)=    z_1\bar{z}_1 + \dots +z_p\bar{z}_p -z_{p+1}\bar{z}_{p+1} -\dots -z_{p+q}\bar{z}_{p+q}.
\end{equation}

We define $G=U(p,q)$
 to be the  unitary group
 that preserves the quadratic form $Q$.  
Let $H^+$ be the stabilizer
 of the vector $e_{1
 }$ and let $H^-$  be the stabilizer of  the vector $e_{p+1}$. 
Then $H^+$ is isomorphic to $U(1,0)U(p-1,q)$ and $H^- $ is isomorphic to $U(p,q-1) U(0,1)$. 
We consider the unitary symmetric spaces
 \[X^+= G/H^+ =U(p,q)/U(1,0)U(p-1,q) \] and \[X^-= G/H^-=U(p,q)/U(0,1)U(p,q-1).\] 

%


\subsection{ } \label{more about Y}
Let $G' =U(p-1,q) $ be the subgroup which fixes the vector $e_{p}$. We consider the  symmetric spaces \[Y^+=G'/G'\cap H^+ = U(p-1,q)/U(1,0)U(p-2,q)\] and \[Y^-=G'/ G' \cap H^- =U(p-1,q)/ U(0,1)U(p-1,q-1)\] 

Consider now the stabilizer $\tilde{G}'\subset G$ of the line through $e_{p}$; then
\[\tilde{G}' =U(p-1,q) U(1,0)\]
It is the fix point set of a involution of G. Furthermore
\[\tilde{H'}^+ = H^+ \cap \tilde{G}' =U(1,0)U(1,0) U(p-2,q)\] and \[\tilde{H'}^- = H^- \cap \tilde{G}' =U(1,0)U(p-1,q-1) U(0,1).\]
Thus we observe 
\begin{eqnarray*}
 \tilde{G}'/ \tilde{H'}^+ &= &U(p-1,q)U(1,0)/U(1,0)U(p-2,q)U(1,0) \\ & = & U(p-1,q)/U(1,0)U(p-2,q) \\ &=& G'/H^+ \cap G' =Y^+
 \end{eqnarray*}
 and 
 \begin{eqnarray*}
 \tilde{G}'/\tilde{H'}^-  & = & U(p-1,q)U(1,0)/U(1,0)U(p-1,q-1)U(0,1)  \\ & = &
  U(p-1,q)/U(p-1,q-1)U(0,1)\\ 
  &= & G'/H^- \cap G' =Y^- 
 \end{eqnarray*}

\medskip
\subsection{}The maximal compact subgroups of $G$, $G'$ and $H^+$ $H^-$ are denoted by  $K, $ $K'$  respectively $K_H^+$ , $K^-_H$. The Lie algebras of the groups are denoted by the corresponding lower case Gothic letters. We use the subscript $_{\bC}$ to denote their complexification.

The representation of the larger group G are always denoted by capital $\Pi$ and those of the smaller group G' by $\pi $.

\subsection{}
{\em As mentioned in the introduction to avoid considering special cases we make in this article the following

{\bf Assumption ${\mathcal O }$:}\label{mathcal O}
  
\begin{center}
              $p \ge 3$ and $q \ge 3$.  
\end{center}
}


\bigskip
\section{The symmetric spaces $X^+$ and $X^-$ and  their discrete series representations}
In  2017 Y.~Sakellaridis and A.~Venkatesh published their book about spherical varieties and periods\cite{SV} which includes several  conjectures of about spherical spaces and the representations in their discrete $L^2 $--spectrum.   The symmetric spaces $X^+$ and $X^-$  are spherical spaces and in this case the conjectures have been verified by  C.~Moeglin and D.~Renard \cite{MR}. In this section we recall and summarize some of their results.

\subsection{Symmetric spaces, their Levi subgroups and $L$-groups }
In \cite{SV} Y.~Sakellaridis and A.~Venkatesh introduced for a spherical space  $X=\hat{G}/\hat{H} $ the concept of a complex  Levi subalgebra of the  complex Lie algebra of $\hat{\bg}_{\bC} $. The symmetric spaces $X^+$ and $X^-$ are spherical spaces and both symmetric spaces have the same complex  Levi algebra $\bl_X=gl(1) ^2\oplus gl(p+q-2,\bC)$  \cite{MR}.

 Y.~Sakellaridis and A.~Venkatesh also defined for a spherical space  an $L$-group.   Following the arguments of C.~Moeglin, D.~Renard, Knop and  Schalke, $ X^+$ and $X^-$ have the same $L$-group which we denote by   $ \check{G}_X  $. It is  $ Sp(2,\bC) \times W_{\bR}$ (page 11 in \cite{MR} ) (The $p+q$  odd is problematic.  C.~Moeglin and D.~Renard  use in their proof of the conjecture of Sakellarides-Venkatesh the group  $ \check{G}_X =SO(2,\bC) \rtimes W_{\bR}$. See remark 8.3 in \cite{MR}. ) All other symmetric spaces for unitary groups have a different $L$-group.



\subsection{ Discrete series representations of the symmetric spaces $X^+$ and $X^-$.}  \label{sec:discrete spec}
\begin{definition} We say that an irreducible representation $\Pi$ of U(p,q)  is a discrete series representation of a symmetric space $X$ for $G$ or equivalently is in the discrete spectrum of a symmetric space $X$ if 
\[Hom_G(\Pi,L^2(X) )\not = 0 .\]
\end{definition}

The work of T.~Oshima and T.~Matsuki \cite{OM} shows that both symmetric spaces $X^+$ and $ X^-$ have discrete series representations. Their Harish Chandra modules are cohomologically induced from $\theta$-stable parabolic subalgebras $\bq^+$ ,$\bq^-$.  
 C.~Moeglin and D.~Renard show that their Levi subgroups are \[ L^+= U(1,0)U(1,0) U(p-2,q),\] respectively 
  \[ L^-= U(0,1)U(0,1) U(p,q-2).\] 
 Assume we have (for the standard torus) a positive root system with only one non compact  {\it simple } root, i.e the simple compact roots are $e_i -e_{i+1}$ for $i= 1,\dots ,p-1$ and $i=q, \dots q-1$ and the simple non compact root is $e_{p}-e_{p+1}$. The representations $\Pi^+$ respectively $\Pi^-$ are parametrized by the characters
(of the Levi subgroups)  \[\lambda ^+= (a_1,-a_1, 0,\dots 0)\] with respectively \[\lambda^- =
(0,\dots ,0, a_2,-a_2)\] with $a_i,$ integers if $n$ is  odd and $a_i$ 1/2 integers if $n$ is even. Furthermore  $a_i \geq \frac{p+q-1}{2}$ \cite{MR}. We write $\Pi^+_{a_1}$ and $\Pi^-_{a_2}$, but if there is no confusion possible we drop the subscript. 
The  minimal $K$-type of $\Pi^+_{a_1}$ is  trivial on $U(q)$ and on the other factor the highest weight is 
\begin{eqnarray*} 
\lambda^+_K &= &(a_1-\frac{p+q-1}{2} +q, 0,\dots,0, -a_1+\frac{p+q-1}{2} -q ) 
\end{eqnarray*}
The minimal $K$-type of $\Pi^-_{a_2}$ is trivial on
$                                  U(p)$ and on the other factor the highest weight is 
\begin{eqnarray*} 
\lambda^-_K &= &(a_2-\frac{p+q-1}{2} +q, 0,\dots,0, -a_2+\frac{p+q-1}{2} -q ) 
\end{eqnarray*}
\cite{KV} To simplify the formulas for the highest weights of the K-types we introduce the notation ${a_0 }= a- \frac{p+q-1}{2}$. 

The representations have the same infinitesimal character if $a_1 =a_2 =a$. In this case we use the notation $\Pi^+_a$ and $\Pi^-_a$.
 

A discrete series representation $\Pi $ of a symmetric space $X \in \{X^+ , \ X^-\} $ has multiplicity \[\mbox{dim Hom}_G(\Pi,L^2(X) )= 1.\]
{ See for example the survey article by E.van den Ban \cite{vB} }

\medskip
\subsection{Arthur packets } An Arthur parameter is a homomorphism $ \phi : W_\bR \times SL(2,\bC) \rightarrow ^L U(n) $ , so it can be considered as a $n$ (here $n = p+q$) dimensional representation of  $W_\bR \times SL(2,\bC).$ 
It is of ``bonne parit\'{e}" if the condition in the first paragraph of 8.3  in \cite{MR} or the previous article  on Arthur packets of unitary groups by C.~Moeglin and D.~Renard \cite{MR2} are satisfied.

D.~Renard and C.~Moeglin verified 
 that the representations $\Pi_{a}^+, \Pi_{a}^-$ in the discrete spectrum of the symmetric spaces $X^+$, respectively $ X^-$, are both in the same Arthur packet  $\mathbb{A}(\phi )$, and that the Arthur parameter is of {\em bonne parit\'{e}}. 
  
  
See more in \cite{MR}

\bigskip

\medskip
\subsection{ The group $ \mathcal{A} (\phi )$ and its characters}
\label{sec: epsilon}
 In \cite{ABV} J.~Adams, D.~Barbasch, and D.~Vogan  defined for an Arthur parameter $\phi $ a group 
\[ \mathcal{A} (\phi ) = \mbox{normalizer } (\phi)/\mbox{centralizer} (\phi ) \]
For a fixed Arthur parameter $\phi $ the representations  in the Arthur packet  of  $\phi $, are parametrized by a characters of  ${\mathcal A}(\phi )$. 
This parametrization depends on the choice of a normalization \cite{ABV}.
In our case
\[  \mathcal{A} (\phi )  =\bZ_2 \times \bZ_2 \]
with generators $E_1, E_2$.
After choosing  a normalization
  C.~Moeglin and and D.~Renard assign  to the  representation in $\{\Pi_a^+,\Pi_a^-\}$ the characters $\epsilon_1, \epsilon_2 $  of ${\mathcal
  A}(\phi)$ 
 \[\epsilon_1(E_1) = 1 ,  \epsilon_1(E_2) = -1 \mbox{ is assigned to } \Pi^+_a\]
 \[ \epsilon_2(E_1) = -1, \epsilon_2(E_2) = 1 \mbox{ is assigned to } \Pi^-_a \]
 and so  abusing the notation we will write sometimes  $\Pi_a(\epsilon_1) $ instead of $\Pi^+_a$ and $\Pi_a(\epsilon_2)$ instead of $\Pi^-$.

Observe that there are  representations in the Arthur packet ${\mathcal A}(\phi)$ which are not in the discrete spectrum of a symmetric space with and they correspond to the remaining 2 characters of $A(\phi)$ .


\subsection{Other descriptions of the symmetric spaces $X^+$ and $X^-$ and their discrete spectrum.} \label{sec: u(n) to su(n)}
We have the exact sequence of groups
\[Z_{p+q}  \rightarrow U(1) \times SU(p,q) \rightarrow U(p,q) \] 
where $Z_{p+q} $ is the set  of the (p+q)'th roots of unity. We identify a root
of unity 
$\zeta \in U_{p+q}$   with the pair \[\{ (\zeta , \zeta I_{p+q}) \subset U(1) \times SU(p,q)\]
in the center of $U(1) \times SU(p,q)$.  $U_{p+q}$ as well as the first factor U(1) act trivially  on $X^+$ and $X^-$  and thus  both  $X^+$ and $X^-$ are  also symmetric spaces for SU(p,q) and have an action of  the covering group $U(1)\times SU(p,q) $.

The restriction of a representation $\chi$ of $U(1)$ to the  group $Z_{p+q}$ of $(p+q)$th roots of unity defines a character $\chi_1$, and an irreducible  representation $\tau $ of $SU(p,q)$  restricted  to the center of $SU(p,q)$ defines a character $\tau_1 $. Then the representation $\chi \times \tau $ of $U(1)\times SU(p,q) $ defines a representation of $U(p,q)$ if $\tau_1 = \bar{\chi}_1$. Here   $\bar{\chi} $ denotes complex conjugate character.  This allows us to  lift a representation of $U(p,q)$ to a representation of $ U(1) \times SU(p,q)$ via this choice of a character. 

Conversely suppose that now that $\check{\Pi}^+$ is a representation in the discrete spectrum of 
$L^2(SU(p,q)/S(U(1)U(p-1,q))$. To show that it extends to a representation $\Pi$ of 
$U(p,q)$ it suffices to show that the center $Z_{p+q}$ of $SU(p,q)$ acts trivially on $\check{\Pi}^+$. Since the center $Z_{p+q}$ acts as a scalar on the irreducible representation $\check{\Pi}^+$ and the center is contained in the maximal compact subgroup it suffices to show that it acts trivially on the highest weight of the minimal $K$-type of $\check{\Pi}^+$. Thus it follows from the formulas in 
Section VI. for the highest weight of the minimal $K$-types of $\Pi^+$ and $\Pi^- $.

\subsection{Discrete series representations of the symmetric spaces $Y^+$ and $Y^-$}\label{pi times trivial}

The discussion of the symmetric spaces  $Y^+$ and $Y^-$ and their  discrete series representations is exactly analogous to the discussion for $X^+$ and $ X^-$. We denote the discrete series representations by the letters $\pi^+$ and $\pi^-$, the Arthur parameter by $\phi'$ and the characters $\epsilon'_i$.

{
Recall that $Y^+$ and $Y^-$ can be considered as  symmetric spaces  of $U(p-1,q-1)U(1,0)$. Hence we will sometimes  identify $\pi^+$ and $\pi^-$ with $\pi^+ \times \text{trivial}$ and $\pi^- \times \text{trivial}$ \ref{more about Y}.

}


\bigskip

\section{Restriction of representations and relative branching laws. }
We introduce in this section   {\em relative branching} and {\em strong relative branching} for representations in the discrete spectrum of symmetric spaces and discuss some examples.

\subsection{Branching laws}

\subsubsection{Unitary branching} Let $\Pi$ be a unitary representation of G. We say that  a unitary representation $\pi$ of a subgroup $G'$ is a subrepresentation of $\Pi_{|G'} $ if there is a nontrivial $G'$-equivariant bounded operator so that  
                      \[ \mbox{Hom}_{G'} (\pi, \Pi_{|G'})\not = 0 \]
i.e if $\pi$ is a direct summand of the restriction of $ \Pi $ to G'. {\it In this article branching will always refer to unitary branching and a {\em unitary branching law} is a rule describing the pairs of irreducible unitary representations $\Pi,\pi$ so that $\pi $ is a direct summand of $\Pi_{|G'}$}.

\subsubsection{Relative branching laws}
Consider   groups $G' \subset G$  with symmetric spaces $Y,X$ with $Y \subset X$ . Let $\Pi $ be a unitary representation in the discrete spectrum of X and $\pi$ a unitary representation in the discrete spectrum of Y. We restrict $\Pi$ to $G'$ and again consider  the G' -equivariant bounded operators 
$\mbox{Hom}_{G'}( \Pi_{|G'} , \pi).$
 A {\em relative  branching law } for representations in the discrete spectrum of the symmetric spaces $X,Y$ with $Y \subset X$ is rule that describes the pairs $\Pi, \pi $ so that 
\[Hom_{G'}( \Pi_{|G'} , \pi) \not = 0.\]

The {\em strong relative branching law.} \label{str_rel_branching} We consider  $\Pi$  a subrepresentation $\Pi$ of $L^2(X) $ and $\pi$ as subrepresentation of $L^2(Y)$ satisfying the relative branching law. The restriction operator
\[ Res_{X\rightarrow Y}: C^{\infty}(X) \rightarrow C^\infty(Y) \] is a symmetry breaking operator. 
Note that at this point the restriction is only defined on sufficiently smooth vectors of $\Pi$, see \cite{OS}.
We could for example work with $K$-finite vectors.
We say that $\pi$  is a strong relative restriction of $\Pi$ if $\Pi,\pi$ is a relative restriction and 
 \[Res_{X\rightarrow Y}( C^{\infty}(X)\cap \pi'  ) \not = 0.\]
The strong relative restriction law is the set of pairs $\Pi,\pi$ 
so  $\Pi,\pi $ a strong relative restriction of $\Pi$.

\medskip
\subsection{ Relative branching for compact unitary symmetric spaces}

The inclusion of the compact symmetric space  Y into the compact symmetric space X  induces a push forward of distributions $D(Y)$ into $D(X)$. An interesting distribution from the representation theoretic perspective is integration of  a function $F(x)$ in $C^\infty(X)$ (suffciently regular in order to make the
integral converge) and 
a function $f(y)$ in $C^\infty (Y) $ which transform according to a irreducible finite dimensional   representation of $G'$ on $C^\infty(Y). $  
 \[ F  \rightarrow \int_Y (f(y)F(y) ) dy \]
 is often referred to as a period integral.
 
Using the example of the compact group $U(n)$ we discuss in this subsection the differences of branching, relative branching, and strong relative branching using period integrals.  We will also consider later other examples  of relative branching for representations for non compact low dimensional symmetric spaces.

For the group  $G=U(n)$ consider the symmetric space (complex projective space) \[X = U(n)/U(1)U(n-1).\] The representations of $G$ are parametrized by their highest weight an $n$-tuple $ (a_1, \dots ,a_n) $ of integers. It defines a representation of $SU(n)$ if the $a_i $ sum up to zero. In particular a representation has an $U(1)U(n-1)$ fixed vector if the highest weight is 
$(a, 0, \dots,0 , - a).$ We restrict the representations to a subgroup $U(n-1)$ which fixes the vector $e_n$. Then $U(n-1) \cap U(1)U(n-1) = U(1)U(n-2) $ and so get a symmetric space \[Y=U(n-1)/U(1)U(n-2) \subset X.\]  The classical branching rule for the restriction to a subgroup $U(n-1)$ shows that all representations of $U(n-1)$ with highest weight $ (b,0,\dots ,c)$ so that $0\leq b\leq a$ and $-a \leq c\leq 0$ are nonzero in the restriction. The representations have a nonzero $U(-1)U(n-2)$ invariant vector iff $c = -b.$ So the relative branching  is expressed in highest weights.

\begin{lemma} (Relative branching law )\\
Let $\Pi$  be a representation in the discrete spectrum of X with highest weight (a,\dots -a).
and $\pi $   be a representation in the discrete spectrum  of Y with highest weight (b, 0 \dots, 0, -b). Then 
\[Hom_{G'}( \Pi_{|G'} , \pi) \not = 0\] 
iff $a \ge b \ge 0$. In this case its dimension is one.
\end{lemma}

\medskip
The strong relative branching law is proved using period integrals and Jacobi polynomials    in
VI.0.2. It states that in this case the strong relative branching law and the relative branching law are identical. Indeed, what we shall see there is that the strong relative restriction law is obtained
simply be restricting a spherical function of $G=U(n)$ to the subgroup $G'=U(n-1)$ and integrating 
over the subgroup against a spherical function
for the subgroup.


\subsection{Relative branching for a non compact unitary symmetric space}

\subsubsection{Discrete series representations}
A familiar case of relative branching is the restriction of discrete series representations of $U(p,q)$ to $U(p-1,q).$ Here we consider the discrete series representations of $U(p,q)$ as discrete series representations of the symmetric space $ X= U(p,q) \times U(p,q)/U(p,q) $ where we embed $U(p,q)$ diagonally into into the product. We define  the symmetric space $Y$ similarly for $U(p-1,q)$  and consider it as an orbit of 
$U(p-1,q)$ on $X.$  Then the discrete spectrum of the unitary restriction of $\Pi $  to $U(p-1,q)$ is isomorphic to its relative discrete spectrum. To see this in general, consider
$$X = G \times G/\text{diagonal }G. $$ 
 Then 
$$ L^2(X) = \sum \Pi \otimes \Pi^* \oplus \int \Pi \otimes \Pi^*$$
where we symbolically write the discrete part of the spectrum as a sum and the continuous part of the
spectrum as an integral. Now the restriction to a reductive $G' \subset G$ will in a similar way
decompose as (for every $\Pi$)
$$\Pi_{|_{G'} }= \sum \pi_i \oplus \int \pi_j$$
(there could be multiplicities here ) 
and hence the tensor product decomposes as
$$(\Pi \otimes \Pi^*)|_{G'} = (\sum \pi_i \oplus \int \pi_j) \otimes (\sum \pi_i \oplus \int \pi_j)^*.$$
Now we see that the spherical part of (i.e. the relative) branching law amounts to taking
the diagonal parts here, and the relative discrete spectrum  becomes exactly, for a fixed
discrete series $\Pi$, the sum
$$\sum \pi_i \otimes {\pi_i}^*$$
which as a discrete part of a branching law is isomorphic to the one for the restriction of the discrete series representation $\Pi$ restricted to $G'$.
Note here the general principle that the continuous spectrum cannot contribute to the discrete
spectrum in this branching -  and so the relative branching law is equivalent to the unitary branching law. The unitary branching law for the  restriction of discrete series representations of $U(p,q)$ to $U(p-1,q)$ was proved by H.He in 2015.\cite{H}

Note that in the case of representations of a compact group $G=U(p)$ and $G'=U(p-1)$  $$X = G \times G/\text{diagonal }G $$ and $$L^2(X) = \sum \Pi \otimes \Pi^* $$ where we sum over all irreducible finite dimensional representations and the character of an irreducible is a function on G.  The analogue
of this period integral considered in the previous example can be written as an identity  of the characters $\chi$ of the representations
$$I(\chi_{\Pi} ,\chi_{\pi}) = \int_{G'} \chi_{\Pi}(x) \chi_{\pi_i}(x^{-1})dx.$$ 
Hence this counts the multiplicity of $\pi_i$ in $\Pi_{|G'}$.

\medskip
\subsubsection{Discrete series representations of other unitary symmetric spaces}
The restriction of a discrete series representation $\Pi^- $ of \[X^-=U(p,q)/U(1)U(p,q-1)\] to \[Y^-=U(p-1,q)/U(1)U(p-1,q-1).\]
for $p>2$ is discussed in section \ref{branchingPi-} and we will see that the relative branching law of  restriction of $\Pi$ is also  equal to the strong relative branching law.

\medskip
\subsubsection{ Whittaker models} Another interesting case is the restriction of  the discrete series  representations of $U(n,n)$ with a Whittaker model to $U(n-1,n)$ respectively $U(n,n-1)$. In this case using He's algorithm we can see that there are representations with a Whittaker model in the discrete spectrum of the restriction of the discrete series representations with a Whittaker model, The discrete spectrum contains also representations without Whittaker model, but never the less we can conclude that the relative branching is non-trivial . 

\medskip

\subsection{Period integrals} \label{ssec:period}
Recall that for irreducible representations $\Pi \subset L^2(G/H)$  and $\pi \subset L^2(G'/H')$ in the discrete spectrum we have 
\[ \mbox{Hom}_{G'}( \Pi_{|G'} , \pi) = \mbox{Hom}_{G'}( \Pi_{|G'} \otimes \pi^{\vee },\bC) \]
where $\pi^{\vee }$ is the contragredient representation of $\pi $. 
 By \cite{AV} $\pi^{\vee} $ is isomorphic to $\pi$.

 We  define a $G'$-invariant linear functional on the $C^\infty $ functions in  $\Pi \subset L^2(X)$ by a {\em  period integral }
 \[ \int_{G'/H'} f(x') f'(x') dx' .        \]
 Here $f$ is a $C^\infty$-function in $\Pi \subset L^2(X)$ and $f'$ is a $C^\infty$-function in $\pi \subset L^2(Y)$. Again, the technical issue here is when this integral makes sense for the appropriate
pairs of functions in the representation spaces, typically $K$-finite resp. $K'$ finite functions,
or smooth vectors, see \cite{OS}. In section \ref{branchingPi+} we will consider the convergence of this integral for pairs $(\Pi^+, \pi^+)$.
 We observe that 
 this linear functional depends only on the restriction of the function $f$ to $G'/H'$. 
 We  summarize this discussion.
 
 \medskip
 {\em Supose that $f$ is a $C^\infty$-K finite function of a representation $\Pi$ in the discrete spectrum of $L^2(X)$ and $f'$ is a $C^\infty$ K' finite function of $\pi $ in the discrete spectrum of $L^2(Y)$. Suppose that a period integral  defines a linear G'-invariant linear functional in 
 \[\mbox{Hom}_{G'}( \Pi_{|G'} \otimes \pi,\bC) \]
 If it is nonzero on such a pair $f,f'$ of $C^\infty$--functions in $(\Pi \times \pi \subset L^2(X)  \times L^2(Y)$ then  the pair $\Pi,\pi$ is a strong relative restriction of $\Pi_{|G'}$.}
\ref{str_rel_branching}

 We will use  period integrals in Section \ref{sec:periods} to obtain a strong relative branching law for the restriction of a representation $\Pi^+_a$ to $G'$.
 
 \medskip

 If $X =G/H$ is a compact symmetric space for a compact group $G$ and $Y$ a symmetric space $G'/(H \cap G')$  period integrals can also help to obtain strong relative branching laws. If $G$ is a classical and $X$ a rank one symmetric space then classical formulas for Jacobi polynomials show that strong relative branching and relative branching coincide. We conjecture that this is always true.


\bigskip

\section{The Main Results}
In this section we state results of the relative branching laws for the restriction of G' of the representations $\Pi_a^+$ and $\Pi_b^-$ in the discrete spectrum of the symmetric spaces $X^+$ respectively $X^-$  and we deduce  theorem
\ref{theorem:1} and \ref{theorem:2} in the introduction.  We then proceed to show that the interlacing patterns of the infinitesimal characters of the representations can be used to define characters on the group  ${\mathcal A}(\phi)$ introduced in \ref{sec: epsilon} and then use this to find a different formulation of the relative branching law \ref{theorem: branch characters}. This reformulation of the  relative branching laws in theorem \ref{theorem: branchPi+} and theorem \ref{theorem: branchPi-}  relates our relative branching laws to  the formulas and ideas of Gross-Prasad and Gan-Gross Prasad (\cite{GP}, \cite{GGP}, \cite{MR}).\\

\medskip

\subsection{Relative Branching laws for the restriction of $\Pi^+$ to $G'$} 
 Let $\lambda=(a,-a, 0,\dots ,0) $ be a parameter of a representation $\Pi_a^+$ of $U(p,q)$ 
 and $\lambda' =(b,-b,0,\dots ,0)$ a parameter for a representation $\pi_b^+$ of $G'=U(p-1,q)$. Recall that $2a$ is an odd integer if $p+q$ is even and an even integer if $p+q$ is odd; $2b$ is an even integer if $p+q-1$ is odd and an odd integer if is even. We denote the representations by $\Pi^+_a$ or sometimes  $\Pi_\lambda(\epsilon_1)$ respectively $\pi_b^+ $ or $\pi_{\lambda'}(\epsilon _ 1)$. Recall assumption $\mathcal O$ \ref{mathcal O}
 
 \begin{theorem}\label{theorem: branchPi+}
  Assume that $G$ and $G'$ satisfy assumption $\mathcal O$ and that $a \geq \frac{p+q-1}{2} $ and $ b \geq \frac{p+q-2}{2} $. Then
 \[ \mbox{dim Hom}_{G'} (\Pi^+ _a , \pi_b  ^+) =1 \]
 if $a> b \geq \frac{p+q-2}{2} $ and zero otherwise. Furthermore
   \[ \mbox{dim Hom}_{G'} (\Pi^+ _a , \pi_b  ^-) =0 \]
for all $b  \geq \frac{p+q-2}{2}$.
 \end{theorem}
 
 
 \medskip
 
 \subsection{Relative Branching laws for the restriction of $\Pi^-$ to $G'$}
  Let $\lambda=( 0,\dots ,0, a, -a) $ be the parameter of a representation $\Pi^-$ of $U(p,q)$ 
and $\lambda' =(b,-b,0,\dots 0)$ the parameter for a representation $\pi^-$ of $G'=U(p,q-1)$. We denote the representations by $\Pi^-_a$ or sometimes  $\Pi_\lambda(\epsilon_2)$ respectively $\pi^-_{b}$ or
 $\pi_{\lambda'}(\epsilon_2) $.

  \begin{theorem} \label{theorem: branchPi-} 
 Assume that $G$ and $G'$ satisfy assumption $ {\mathcal O}$ and that $a \geq \frac{p+q-1}{2} $ and $ b \geq \frac{p+q-2}{2} $. Then
 \[ \mbox{dim Hom}_{G'} (\Pi^- _a , \pi_b  ^+) =0 \]
 for all $b\geq \frac{p+q-2}{2}$. Furthermore 
  \[ \mbox{dim Hom}_{G'} (\Pi^- _a , \pi_b  ^-) =1 \]
iff $b>a \geq \frac{p+q-1}{2} $ and zero otherwise. 
 \end{theorem}
 \
 
Combining relative branching laws for $\Pi^+ $ and  $\Pi^-$ proves the  theorem \ref{theorem:1} in the  introduction.

Theorem \ref{theorem: branchPi+} and \ref{theorem: branchPi-} are proved in sections VI. and VII.

 \subsection{Interlacing patterns and characters of the group ${\mathcal A}(\phi)$}
 The Gross Prasad Tan conjectures for discrete series  as proved by 
 R. Beuzart-Plessis are phrased using a character on the group ${\mathcal A}(\phi) =(Z_2)^{p+q}$. We are using these ideas and will start by relating the characters  to related interlacing properties the parameters of the representations.

 \begin{definition} \label{def: RB} We say that   a pair a,b  of half integers has the {\em property RB } if 
 \begin{enumerate} 
 \item  2a = 0 mod 2  and 2b = 1 mod 2  if p+q is odd\\
 and 2a=1 mod 2 and 2b = 0 mod 2  if p+q is even
 \item $a-\frac{p+q-1}{2} \in \bN$ and $b-\frac{p+q-2}{2} \in \bN$
 \end{enumerate}
 \end{definition}
Ordering the  numbers $a,-a,  b,-b$
 in decreasing order we have  the interlacing patterns 
 \begin{eqnarray*}
     P_1(a,b) &:   &  (a, b,-b, -a)   \\
      P_2(a,b) &:   & (b,a , -a, -b)
\end{eqnarray*}

For an interlacing pattern we define a characters $\tilde{\epsilon}$,  $\tilde{\epsilon}'$ of ${\mathcal A}(\phi )=\bZ_2 \times \bZ_2$ with generators $E_1, E_2$ as follow. Let $a(i,>)$ be the number of $b_j $ in the interlacing pattern larger than $a_i$ and $b(j,>)$ the number of $a_i $ larger than $b_j$ in the pattern. Following the ideas of B.~Gross and D.~Prasad \cite{GP} we define characters on ${\mathcal A}(\phi)$ as follows

\begin{definition} Suppose that $a,b$ have the property RB. Define 
\begin{eqnarray*}
\check{\epsilon}(E_i) &= &(-1)^{i+1+ a(i,>)}  \\
\check{\epsilon}'(E_i) &= & (-1)^{i+ b(j,>)}  
\end{eqnarray*}
\end{definition}

So  if $a, b$ have the interlacing pattern 
$P_1$ then  \[\check{\epsilon}  =\check{\epsilon}' = \epsilon_1\] (see  \ref{sec: epsilon} for a definition of the characters $\epsilon_1,\epsilon_2$). 
If $a, b$ have the interlacing pattern $P_2$ then 
 \[\check{\epsilon}  =\check{\epsilon}'= \epsilon_2 \]
 
 So we have a correspondence between pairs (a,b)  with the interlacing patterns $P_1$ $P_2$ and pairs of representations 
 \[P_1(a,b) <-->(\Pi_a(\epsilon_1), \pi_b(\epsilon'_1) ) = (\Pi_a^+,\pi_b^+).\]
 \[P_2(a,b) <--> (\Pi_a(\epsilon_2), \pi_b(\epsilon'_2) ) = (\Pi_a^-,\pi_b^-).\]
 

\subsection{Interlacing patterns and relative branching}
We  reformulate  Theorem \ref{theorem:1} and Theorem \ref{theorem:2}
using the considerations in \ref{def: RB}
 as follows

\begin{theorem} \label{theorem: branch characters}

Let  $a, b$ be positive numbers with property RB.
\begin{itemize}
\item Suppose that $a, b$  satisfy interlacing pattern $P_1$  and thus correspond to  the pair $\Pi_a(\epsilon_1), \pi_b(\epsilon'_1) $. The representations $\Pi_a(\epsilon_1)$   and  representation $\pi_b(\epsilon'_1)$    satisfy  \[ \mbox{dim Hom}_{G'} (\Pi _a(\epsilon_1) , \pi_b(\epsilon'_1) ) \not = 0 \] 
\item  Suppose that $a, b$ satisfy interlacing pattern $P_2$ and thus correspond to the pair  $\Pi_a(\epsilon_2), \pi_b(\epsilon'_2) $. The representations $\Pi_a(\epsilon_2)$   and   $\pi_b(\epsilon'_2)$    satisfy  \[ \mbox{dim Hom}_{G'} (\Pi_a(\epsilon_2) , \pi_b(\epsilon'_2)) \not = 0 \] 
\item If $i \not= j$ 
\[ \mbox{dim Hom}_{G'} (\Pi_a(\epsilon_i) , \pi_b(\epsilon'_j))  = 0 \] 

\end{itemize}
\end{theorem}

An equivalent direct formulation avoiding the characters of ${\mathcal A}(\phi)$

\begin{cor} \label{cor: simple} 
Let  $a,\ b$ be positive numbers with the property RB
\begin{itemize}
\item $a, b$  satisfy interlacing pattern $P_1$ if and only if
 \[ \mbox{dim Hom}_{G'} (\Pi _a^+ , \pi_b^+)  \not = 0 \] 
\item $a, b$  satisfy interlacing pattern $P_2$ if and only if
 \[ \mbox{dim Hom}_{G'} (\Pi _a^- , \pi_b^-)  \not = 0 \]  
\end{itemize}
\end{cor}

\begin{remark}
Although the formulation in  corollary \ref{cor: simple} is simpler, the formulation in \ref{theorem: branch characters} suggests generalizations to other unitary symmetric spaces and other groups. A similar approach can be used to reformulate the results in the article ``A hidden symmetry of a branching law". \cite{KS}.
Note that in our case the dimensions above are exactly one when non-zero.
\end{remark}

\medskip

\subsection{A multiplicity one result}

Following the ideas of B.~Gross and D.~Prasad in \cite{GP} we  reformulate  these results.

Assume again that $
a,b$ satisfy \ref{def: RB}
Define 
\[A_S (a) = \Pi^+_a \oplus \Pi^-_a \]

Similarly we define    
\[B_S'(b)= \pi^+_{b} \oplus \pi^-_{b}\]
  Following the ideas  in \cite{GP}  we consider for $\pi_{b} \in \{ \pi_b^+,\pi_b^-\}$
\[ \mbox{dim Hom}_{G'} (A_S(a),  \pi_b)\]
and for $\Pi _a \in \{\Pi_a^+, \Pi_a^- \} $
\[ \mbox{dim Hom}_{G'}(\Pi_a, B_S(b)  )\]
and deduce

\medskip
\begin{cor}\label{cor: mult1}
Assume that  $p,q > 3$ and $p \not = q$, Assume again that a,b have property RB
\begin{enumerate}
\item There exists exactly one representation $\pi_b \in \{ \pi_b^+,\pi_{b}^-\}$ so that 
 \[\mbox{dim Hom}_{G'}(A_S(a), \pi_{b}  )\not = 0\]
\item There exists exactly one representation  $\Pi _a \in \{\Pi_a^+, \Pi_a^- \} $ so that 
\[ \mbox{dim Hom}_{G'}(\Pi_a , B_S(b) \not = 0 )\]
\end{enumerate}

\end{cor}

\medskip

{This implies}

\begin{cor}
Assume that  $p,q > 3$ and $p \not = q$ and  that a, b have the property RB . Then
\[\mbox{dim Hom}_{G'}(A_S(a), B_S(b)  )= 1.\]
\end{cor}

\medskip

See \ref{theorem:2} in the introduction.

 
 \begin{remark}
 The example in IV.3 shows that the assumption $p,q > 3 $ is necessary, but the assumption that 
$ p \not = q$ is not essential but very helpful for the exposition.
  \end{remark}
  \bigskip

\section{Relative branching for $\Pi^+$}\label{branchingPi+}

{
In this section we use the ideas in \ref{sec:periods} and construct  explicit non vanishing period integrals on the unitary symmetric spaces  $X^+$ and $Y^+$
thus  obtaining a $G'$-invariant pairing between the discrete series $\Pi_{\lambda}^+$ and $\pi_{\lambda'}^+$ . In the last part 
\ref{sec:exhaustion} we show that there are no further discrete summands in the relative branching
law for the discrete series of the unitary symmetric space in addition to the
ones found using period integrals, thus obtaining a proof of their relative branching law
 in  Theorem \ref{theorem: branchPi+}. Here we use branching in stages along the 
 the sequences
\begin{itemize}
\item $U(p-1,q) \subset O(2p-2, 2q) \times O(2) \subset O(2p, 2q)$
\item $U(p-1,q) \subset U(p,q) \subset O(2p,2q)$
\end{itemize}


\subsection{Background} 
\subsubsection{Affine symmetric spaces}
Recall the
general setup of the affine symmetric space $X = G/H$ and the maximal compact subgroup
$K \subset G$ with the decomposition $G = KBH$ coming from a certain Cartan subspace  $B$ in the
tangent space at the base point for $G/H$. Our case is rank one, i.e. $\dim B$ = 1 and we get the coordinates $K/K_H B H$ where $K_H = U(p-1)U(q)$. 
We also have the extra condition that $G' = K'BH'$, i.e.the two spaces $X^+$ and $Y^+$ have the same hyperbolic variable from $B$.

 \subsubsection{Flensted-Jensen representations and Flensted-Jensen functions}
 
In 1980 M. Flensted-Jensen \cite{FJ} discovered a  family of irreducible unitary representations  in the discrete spectrum $L^2(X)$ of some symmetric spaces for {\em connected  semisimple } Lie groups. Although $U(p,q)$ is not semisimple the considerations in \ref{sec: u(n) to su(n)} show that there is a 1-1 correspondence between the representations in the discrete spectrum of $U(p,q)/U(1)U(p-1,q)$ and $SU(p,q)/S(U(1)U(p-1,q))$ 
 and that we can extend the Flensted--Jensen result to representations in  the discrete spectrum of  $L^2(X^+)$ for $p>1$. Flensted-Jensen's proof is based on an explicit formula for  functions $\psi$ in the lowest $K$-type  of  representations in the discrete spectrum of $L^2(X^+)$  \cite{FJ}. We refer to these functions as Flensted-Jensen functions and the corresponding representations as Flensted-Jensen representations.
The Flensted-Jensen representations are cohomologically induced from a character of a $\theta$--stable parabolic with Levi subgroup $S(U(1)U(1)U(p-2,q))$
and their Langlands parameters were determined by H.~Schlichtkrull \cite{Sch2} \cite{Sch}.





We recall now the explicit construction of M.~Flensted-Jensen functions, i.e. the lowest
$K$-types of the Flensted-Jensen representations  in  the discrete spectrum of $L^2(X^+)$. In this  we   follow here the notation in \cite{FJ} and in particular the notation in the last section of \cite{FJ} where
rank one cases are explained. In particular we shall follow in subsections the notation there in considering the
symmetric space $X = SU(p,q+1)/S(U(p,q) \times U(1))$, i.e. a switch from the
space $ SU(p+1,q)/S(U(p,q) \times U(1)) $  considered above to $X$,  after applying an outer isomorphism $SU(p,q+1) \rightarrow SU(q+1,p)$. This allows us to cite the explicitly the formulas of M.~Flensted-Jensen. We determine in the subsections VI.1 IV.2 using a period integral representations in the discrete spectrum of the restriction to  
$G'= SU(p,q)$ in the discrete spectrum of $Y=G'/S(U(p,q) \times U(1))$.

\medskip
\subsubsection{Period integrals}  \label{sec:periods} 
We  discuss in this section the period integral
$$I(f, f') = \int_{Y} f(y') f'(y')dy'$$
for functions $f(x)\in \Pi^+, f'(y) \in \pi^+$
and determine whether the integral  converges, and whether it is non-zero for the Flensted-Jensen functions $\psi$, $\psi'$. Here
$Y = G'/H'$ and $H' = G' \cap H$. 
The bilinear form $I(\ ,\ )$ is clearly
$G'$-invariant, and hence gives information about the branching of the representation
generated by restricting $f$  from $G$ to $G'$ i.e the {\it relative branching} in the sense of
finding summands in the branching law of the same type, i.e. Flensted-Jensen representations
in the discrete spectrum of $Y$.  For more about period integrals see \ref{sec:periods}.
 
 We assume in  the  considerations and calculations of Flensted Jensen functions and period integrals that 
$$G = SU(p,q+1), q > p > 0, \ H = S(U(p,q) \times U(1))$$
and $H$ fixing the last coordinate. 
Note that now $K=S(U(p) \times U(q+1))$ and 
$$K/K_H = SU(q+1)/S(U(q) \times U(1)) = P^q(\C)$$
the projective space.

\medskip
\subsection{Flensted-Jensen functions and Jacobi polynomials}

We will now find the explicit formulas for the Flensted-Jensen functions $\psi_{\lambda}$ on $X$ corresponding to
the lowest $K$-type which generate the discrete series representation $\Pi^+_{a_0}$ in $L^2(X)$. Note that
by construction the lowest $K$-type in this representation is spherical for $K/K_H$.
Then we shall find the corresponding functions on $Y$.

Recall the  Jacobi polynomials $P_n^{(\alpha, \beta)}(x)$ of degree n which are usually normalized by
 $$P_n^{(\alpha, \beta)}(1) = \frac{\Gamma(n+\alpha+1)}{\Gamma(n+1)\Gamma(\alpha+1)}.$$ 
 
 \underline{Fact} The Flensted-Jensen functions on $X$  are of product type, with $s$ the radial variable,
$$\psi_{\lambda}(k b_s) = (\cosh s) ^{-i\lambda - \rho} P_n^{(q-1,0)}(k)$$
for $k \in K,\, b_s \in B$. Here the parameters are (notation as in \cite{FJ})
$i\lambda = q-p+n \,(n \in 2\Z^+), \rho = p+q, \rho_{\bt} = q,  \rho - 2\rho_{\bt} = p-q,
\mu_{\lambda} = i\lambda + \rho -2\rho_{\bt} = n.$

The Jacobi polynomials $P_n^{(\alpha, \beta)}(x)$ with $\alpha = q-1, 
\beta = 0$ can be considered as spherical functions
 on the compact Riemannian symmetric space
$P^q(\C)$, see \cite{AT} \cite{HS}.                                                                                                                                                                                                                                                                                                                                                                  
An important classical relation in our situation is
$$P_n^{(\gamma, 0)}(x) = \frac{\Gamma(n+1)}{\Gamma(n+\gamma +1)}\sum_{k=0}^n \frac{\Gamma(k+\alpha+1)(2k+\alpha+1)}{\Gamma(k+1)}P_k^{(\alpha, 0)}(x)$$
where $\gamma = \alpha+1 = q-1$. 
This formula is crucial in obtaining  the relative branching
law for spherical representations from $SU(q+1)$ to $SU(q)$ on the complex projective spaces. 

\subsection{Computation of the period integral}
 
We will now compute the period integrals of Flensted-Jensen functions  by considering the
branching on the (maximal compact subvarieties) $P^{q}(\C)$ resp. $P^{q-1}(\C)$.



\subsubsection{The nonvanishing of the period integral} We consider the branching to $G' = SU(p,q)$, fixing the $p+1$ coordinate and have
similar formulas for the Flensted-Jensen functions  on $Y$. The corresponding period
integral is now
\begin{multline} 
I( \psi_\lambda, \psi'_{\lambda'})  = \int_0^\infty \int_{P^{q-1}(\C)} P_n^{(q-1,0)}(y) P_k^{(q-2,0)}(y) \cdot\\
\cdot (\cosh s)^{-2q - n} (\cosh s)^{-2q +2  - k} dy D(s) ds
\end{multline}
where $dy$ is the invariant measure, and the density \[D(s) = (\cosh s)^{2q-1} (\sinh s)^{2p-1}.\]
Here $i\lambda' = q-1-p+k \, (k \in 2 \Z^+)$ is the Flensted-Jensen parameter for $G'$.
The exact value may be calculated in terms of Gamma functions. Note that 
$$A(\alpha,\beta) = \int_0^{\infty} (\sinh t)^{\alpha} (\cosh t)^{-\beta} dt = \frac{1}{2}
B((\alpha - 1)/2, (\beta - \alpha)/2)$$
where the beta function $B(x,y) = \Gamma(x) \Gamma(y)/(\Gamma(x+y)$. So 
\begin{multline} I( \psi_\lambda, \psi'_{\lambda'}) = \\ =\int_0^\infty \int_{P^{q-1}(\C)} P_n^{(q-1,0)}(y) P_k^{(q-2,0)}(y) 
(\cosh s)^{-2q +1  - n - k} (\sinh s)^{2p - 1}dy ds\\
= A(2p - 1, 2q +n + k -1) \int_{P^{q-1}(\C)} P_n^{(q-1,0)}(y) P_k^{(q-2,0)}(y) dy
\end{multline}
and the
non-vanishing is governed by the classical formula above. Note also that the convergence is assured by
$2p-2q-n-k < 0 \,(n \in \N, \, k \in \N).$ And observe finally, that $I$ is non-vanishing
exactly for $0 \leq k \leq n.$ (The degree of $P_n(x)$ is $n$, which here is half the parameter
$n$ in the Flensted Jensen parametrization - see the example below.) 
Here the parameters are $n, k$ and they are exactly the parameters $a, b$ in the previous
discussion.

Thus restricting the spherical function on the 
left hand side to the subgroup and
integrating against a similar spherical function for the subgroup, we get a non-zero period
integral precisely for the representations in the branching law; hence the strong relative branching
is the same as the branching law.
See formula (3.40) in \cite{AF}, which also covers the
relevant relative branching laws for the analogous real and quaternionic case (see the appendix).


\subsubsection{Example:} The Jacobi polynomials corresponding to the projective line $P^1(\C)$
are just Legendre polynomials, and the spherical representations in this case
are the representations of $SU(2)$ with integer highest weight $a$ of dimension
$d = 2a+1$. We have then for the parameter $n$ in the Flensted-Jensen labels
$n = 2a$, and the degree of the polynomial is $a$. In more detail recall the representations
in the model where the space is
$$V_m = span \{z^m, z^{m-1}w, z^{m-2}w^2, \dots, w^m\}$$
i.e. homogeneous polynomials in two variables $z, w$. The action is the
linear action in the variables, and the torus $diag(e^{i\theta}, e^{-i\theta})$
acts with the eigenvalues $\{e^{im\theta}, e^{i(m-2)\theta}, \dots, e^{-im\theta}\}$
and the representation is spherical for $m \in 2\N$. The spherical vector is 
$z^nw^n, \, m = 2n$, and the spherical polynomial is then obtained as the matrix coefficient
of this vector under the rotations by the angle $\theta$, i.e. it is proportional to
$$( z^nw^n, (\cos\theta\sin\theta z^2 + ((\cos\theta)^2 - (\sin\theta)^2)zw - \cos\theta\sin\theta w^2)^n)$$
using the inner product in $V_m$. For example when $n=2$ we get (a constant times)
$3(\cos2\theta)^2 - 1$, which is the Legendre polynomial of degree 2 in the variable
$x = \cos2\theta$; this is consistent with this being the spherical harmonic on the
2-sphere with the polar angle $2\theta$ corresponding to the double covering
$SU(2) \mapsto SO(3)$.

\subsubsection{Conclusion:}
The relative branching of the discrete series $\Pi_{\lambda}$
with lowest $K$-type $\psi_{\lambda}$ for $X$ to the discrete series $\pi_{\lambda'}$
with lowest $K$-type $\psi'_{\lambda'}$ for
$Y$ is by the period integrals containing the parameters (in terms of the degree of the
corresponding Jacobi polynomials $n$ resp. $n'$) $0 \leq n' \leq n$. 
This is the
discrete spectrum in the relative branching theorem, and in terms of the infinitesimal
characters corresponds to the interlacing pattern described in  Theorem V. I.

\subsection{Exhaustion and branching law} \label{sec:exhaustion}
In this section we prove the last key step in the branching law, namely the
exhaustion of the discrete spectrum that we found above via period integrals and show that there are no further discrete summands in the relative branching
law for the discrete series of the unitary symmetric space in addition to the
ones in the conclusion above. We shall combine recent results by T. Kobayashi (2021)
and  Schlichtkrull, Trappa, Vogan (2018) \cite{STV} and do the branching law in stages, based on the
two sequences
\begin{itemize}
\item $U(p-1,q) \subset O(2p-2, 2q) \times O(2) \subset O(2p, 2q)$
\item $U(p-1,q) \subset U(p,q) \subset O(2p,2q)$
\end{itemize}
and with notation as in these references, i.e we consider again the symmetric space $X^+= U(p,q)/U(1)U(p-1,q)$ and restrict to $G'=U(p-1,q)$.

\subsubsection{Geometric considerations} 

We think of the  real hyperboloid (over the real numbers) $H_{p,q}(\R) = O(p,q)/O(p-1,q)$ 
 as
$$H_{p,q}(\R) = \{(x_1, x_2, \dots, x_n) \in \R^n | x_1^2 + \cdots + x_p^2 
- x_{p+1}^2 - \cdots - x_n^2 = 1 \}$$
where $n = p+q$. Similarly we have the complex hyperboloid in complex space
$$H_{p,q}(\C) = \{(z_1, z_2, \dots, z_n) \in \C^n | |z_1|^2 + \cdots + |z_p|^2 
- |z_{p+1}|^2 - \cdots - |z_n|^2 = 1 \}$$
again with $n = p+q$. We identify $H_{2p,2q}(\R)$ and $H_{p,q}(\C)$,
and the second one with $U(p,q)/U(p-1,q)$ (say $U(p-1, q)$ fixing the first coordinate).
If we take the projective complex hyperboloid, i.e. we identify points on the
circles obtained by the action of $U(1)$, then we obtain the symmetric space
\begin{multline}
H_{p,q}(\C)/U(1) = U(p,q)/(U(1) \times U(p-1,q)) \\ = SU(p,q)/S(U(1) \times U(p-1,q)).
\end{multline}

In this way we may identify the discrete spectrum on this symmetric space $X^+$ with 
the $U(1)$ invariant discrete 
spectrum of the real hyperboloid. Note that this $U(1)$ is exactly the center of $U(p,q) $, see \ref{sec: u(n) to su(n)}.


\subsubsection{Relative branching from SO(2p,2q) to U(p,q)}\label{subsec:SO(2p,2q)-U(p,q)}
This is discussed in detail in Theorem 6.1 in \cite{K1}
The last section of \cite{STV} contains also an account of these results.
We summarize the results.
 Here the representations in the discrete in the discrete spectrum  of $H_{2p,2q}(\R)$  are denoted by $U(\ell) = U^{S0(2p,2q)}(\ell)$.
   They are cohomologically induced from a maximal $\theta$-stable parabolic subgroup  with Levi  subgroup $SO(2,0) SO(2p-2,q)$ and parametrized by an integer $\ell > p+q-1$. Here $\ell-p-q+1 $ defines a character of $SO(2,0)$.
 
 
 {
  By Theorem 6.1  in \cite{K2}
  \[ U(\ell)_{|U(p,q)} = \oplus _{x,y}\Pi^+_{(x,y)}.\]
  The representations $\Pi^+_{(x,y)}$ are cohomologically induced 
  from a parabolic subgroup with Levi subgroup $U(1,0)U(1,0)U(p-2,q)$. 
 Here $(x,-y)$ define characters of $U(1)U(1)$ and they satisfy
   $x+y = \ell $ . Then $\Pi^+_{(x,x)} $ is in the discrete spectrum of 
   of the symmetric space $U(p,q)/U(1)U(p-1,q)$ iff  $x=y=\ell/2. $
   
  \subsubsection{ More on  branching from O(2p,2q) to U(p,q)}
 A representations $\chi_s$ of $U(1)$ defines a line bundle on
 $U(p,q)/U(1)U(p-1,q)$0.
T. Kobayashi considers the discrete spectrum of $L^2(U(p,q)/U(1)U(p-1,q),\chi_s)$.
 In \cite{K1}  he  shows that the representations $\Pi^-_{(x,y)}, \, x-y =s$  are in the discrete spectrum of $L^2(U(p,q)/U(1)U(p-1,q),\chi_s)$. Thus  
  \[ U(\ell)_{|U(p,q)} \subset \oplus_{|s|<\ell} L^2(U(p,q)/U(1)U(p-1,q),\chi_s)\]

  \subsubsection{Branching from $SO(2p,2q)$ to $SO(2,0)SO(2p-2,2q)$} \label{subsec:SO(2p,2q) to SO(2,0)SO(2p-2,2q)}
We will use theorem 1.1  in \cite{K}.  We consider the case $p' =2 , q'=0$.  So  we have $\delta =+$. For the trivial representation of SO(2,0) we have $\lambda'= 0$. Furthermore recall that  we have by definition $\ell = \lambda$.  Hence the parameter of the representations in the restriction of $U(\ell)$  are   the sets
\[ \Lambda_{++}(\ell ) \text{ and } \Lambda_{+,-} (\ell) \]
\begin{enumerate}
\item The set $\Lambda_{++}(\ell)$ is finite or empty. If $\lambda'=0$   the representations in the restriction are  discrete series representations of $SO(2p-2,2q)/SO(2p-1,2q)$. We have 
$\ell - \ell" -1 \in 2\bN$.
\item $ \Lambda_{+,-} (\ell) $ is empty if  $\lambda' = 0$.
\end{enumerate}
In this case there is also continuous spectrum in the restriction.

  \subsubsection{Branching from $SO(2p,2q)$ to $SO(2p,2q-2)SO(0,2)$}\label{subsec:SO(2p,2q) to SO(0,2)SO(2p,2q-2)}
 We will use again theorem 1.1  in \cite{K}.  We consider the case $p'=0 , q'=2$.  So  we have $\delta =-$. For the trivial representation of SO(0,2) we have $\lambda"= 0$. Furthermore recall that  we have by definition $\ell = \lambda$.  Hence the parameter of the representations in the restriction of $U(\ell)$  are   the sets
\[ \Lambda_{-+}(\ell )  \]
  This set is infinite so
  that in this case  $U(\ell)_{|O(2p, 2q-2)}$ is a direct sum of irreducible representations. If $\lambda' = 0$ then $\ell " -\ell -1 \in 
  2\bN $.

 \subsubsection{Branching  from U(p,q) to U(1)U(p-1,q)}\label{subsec:U(p,q)-U(1)U(p-1,q)}
 Observe that $U(1)$ commutes with $U(p-1,q)$ , hence  its representation is a character. It also acts by this character on the highest weight vector of the minimal $K$-type. From \cite{STV} theorem ? we deduce 
 \begin{itemize}
 \item $U(1)$ acts trivially on discrete spectrum of  the restriction of  the representation 
$\Pi^+_{(x,x)} $ to  $U(p-1,q)$, 
 \item If $x\not = y$ , $U(1)$ does not act trivially on discrete spectrum of  the restriction of  the representation 
$\Pi^+_{(x,x)} $ to  $U(p-1,q)$, 
\end{itemize}

 \subsection{Branching along the first sequence}
 Combining \ref{subsec:SO(2p,2q) to SO(2,0)SO(2p-2,2q)} and \ref{subsec:SO(2p,2q)-U(p,q)} we conclude that the only representations which appear in the  restriction of $U(\ell)$ to $U(p-1,q)$ along the first sequence and which are in the discrete spectrum of $U(p-1,q)/U(1)U(p-2,q)$ are the representations which we obtained using period integrals.

 \subsection{Branching along the second sequence}
 Combining \ref{subsec:SO(2p,2q)-U(p,q)}, and \ref{subsec:U(p,q)-U(1)U(p-1,q)}
 we deduce that the relative discrete spectrum of the restriction of $U^{SO(2p,2q)} (\ell)$ to 
 $U(p-1,q)$ is the relative discrete spectrum of the representation 
 $\Pi^+_{(x,x)}$.
 
 \medskip
 Combining the  results of the last 2 subsections and the theorem by T. Kobayashi \cite{K3} regarding continuous spectrum, 
 this concludes the proof of branching theorem for the relative restriction of $\Pi^+$.}

We give here a few more details on the branching in stages argument, repeating the
crucial steps:
\subsubsection{Branching along the first sequence} We perform the branching (using \cite{STV}, Section 8 and 11, and \cite{K},
Theorem 1.1) in stages, using that it may be performed in two ways from $O(2p,2q)$  to $U(p-1,q)$,
namely along the first sequence, respectively along the second sequence.

 In the first sequence the
restriction of a discrete series $\Pi^{O}$ for $H_{2p,2q}(\R)$ to $O(2) \times O(2p-2,2q)$ has both a continuous spectrum, and also a discrete
spectrum \cite{K}.  Note that it is a direct sum with the summands  parametrized by characters of $O(2)$. We interested representations of $ U(p,q)$, respectively $U(p-1,q)$ 
with a trivial action of the center   of $U(p,q)$, respectively $U(p-1,q)$ and
representations of $U(p-1,q) \times U(1)$  with trivial action of the  the $U(1)$ \ref{pi times trivial}
in the upper left corner. This implies that we  consider only  the summand which is invariant under the subgroup $SO(2)$ of
 $O(2) O(2p-2,2q)$).  
 
 Claim: The direct summand of $\Pi^O _{|O(2) O(2p-2,2q)}$ which is invariant under $O(2)$ contains  only  finitely many irreducible representations in its discrete spectrum .\\
  To see this we recall
the details of the orthogonal branching from \cite{K} relevant here in our case
$$\tilde{G'} = O(2,0) \times O(2p - 2, 2q) \subset G = O(2p,2q)$$
with the $O(2,0) = O(2)$ acting on the coordinate $z_1$. The discrete spectrum contains
a priori three types of representations, namely for $p' = 2, q' = 0, p'' = 2p - 2, q'' = 2q$
$$(-+) \,\pi_{-,\lambda'}^{p',q'} \otimes \pi_{+, \lambda''}^{p'',q''}$$
which in this case is empty; second case is
$$ (++)\, \pi_{+,\lambda'}^{p',q'} \otimes \pi_{+, \lambda''}^{p'',q''}$$
which is the relevant cased for us (the relative case, where the subrepresentations
are again discrete series for an imbedded real hyperboloid of smaller dimension), here
$\lambda - \lambda' - \lambda'' -1 \in 2\N$. Also the parameters are
$\lambda'' \in \Z + (p'' + q'')/2, \,\lambda'' > 0$, and $\lambda' \in \Z, \, \lambda' \geq 0$.
The parameter $\lambda \in \Z, \, \lambda > 0$. (Note that in \cite{K} the parameter  used is
$\lambda = l +(p+q-2)/2$ for $O(p,q)$ whereas in \cite{STV} is used the parameter $l$.)
The third case is similarly $(+-)$ and involves $\pi_{-,\lambda''}^{p'',q''} = \pi_{+,\lambda''}^{q'',p''}$,
and does not correspond to discrete series for the smaller embedded real hyperboloid.
Note that the representations occuring for the $O(2)$ are spherical harmonics, since in general
$\pi_{+,\lambda}^{p,0} = {\mathcal H}^m(\R^p)$, spherical harmonics of degree $ m = \lambda - \frac{p}{2}
+1$. Hence for the relative discrete branching we get
$$\left(\pi_{+,\lambda}^{2p,2q}|G'\right)_{relative - discrete} = \oplus
 \pi_{+,\lambda'}^{2,0} \otimes \pi_{+, \lambda''}^{2p - 2,2q}$$
sum over $\lambda - \lambda' - \lambda''-1 \in 2\N$ which is a finite sum. The terms with $\lambda' = 0$ constitute the part that we will decompose further under $U(p-1,q)$.

Note that in the
continuous spectrum  in the branching from $G$ to $G'$ there can be no discrete spectrum for $U(p-1,q)$.
To see this, suppose that in the continuous spectrum for the Hilbert space one representation had a restriction
with a discretely occurring representation; then these would consist of vectors in the
Hilbert space, a contradiction, since these only are weakly contained here.
Another way to say this is that the direct integral decomposition of a unitary representation
is with a unique measure class, and hence the discrete part is unique. 

\subsubsection{Branching along the second sequence}
On the other hand,
in the second sequence the first restriction yields the discrete sum over $\pi_{x,y}^{U(p,q)},
\, x+y = l, \, x,y \in \Z$. But since the center of $U(p,q)$ acts by the character of weight
$x-y$ the relative part has $x = y$ so we only get a finite sum.
Finally we obtain the (finite, as it turns out) relative 
discrete spectrum restricting $\pi_{x,x}^{U(p,q)}$
to $U(p-1,q)$ consisting of the sum of $\pi_{a,b}^{U(p-1,q)},\,
a = b \leq x$ in agreement with the exhaustion that we wanted.
\qed

{\bf Remark 1}. Using the same results from \cite{K} and \cite{STV} for the opposite case
of restricting the same representation $\pi_{+,\lambda}^{2p,2q}$
to $G' = O(2p, 2q - 2) \times O(0,2)$, we would find a purely discrete spectrum
(this an admissible restriction) and we would find the full discrete spectrum
of $U(p,q-1)$ inside $U(p,q)$ in agreement with what we found earlier.

\medskip
The non vanishing of the period integral and the above results about  about branching in stages imply 

\begin{theorem} Theorem \ref{theorem: branchPi+} is a strong relative branching law.
\end{theorem}

{\bf Remark:}
As mentioned in  Appendix 2, the quaternionic case could have
been treated in a similar way via $G' = O(4,0) \times O(4p - 4, 4q) \subset O(4p,4q)$.

\bigskip




 \section{Relative branching for $\Pi^-$} \label{branchingPi-}
 
 We prove in this section the  main theorem about the relative branching of $(\Pi ^-,\pi^-)$ for an irreducible representation $\pi^- \subset L^2(G'/H')$ with regular infinitesimal character. The proof in this section uses also restriction in stages and the same ideas as in the previous section. It is also possible to prove it using pseudo dual pairs and the ideas in \cite{S}.
  This theorem was also proved using different ideas earlier by Y. Oshima. We thank him for informing us about this unpublished result \cite{O}. 

\bigskip

\subsection{The representation $\Pi^- $. }
The representation $\Pi ^- $ is cohomologically induced from a  $\theta$-stable parabolic subalgebra with Levi subgroup $U(p,q-2)U(0,1)U(0,1)$ and has parameter  $(0,\dots 0,a,-a)$ which is trivial on $U(p,q-2)$ .
 We assume that it is in the "good range", i.e that $a -\frac{p+q -1}{2} \geq 0$.

The infinitesimal character of $\Pi^-$ is 
up to  conjugation with an element in the Weyl group 
\[(a, \frac{p+q-2}{2} , \dots,- \frac{p+q-2}{2}, -a)\] and it is the same as the infinitesimal character of $\Pi^+$.  
The minimal $K$-type of $\Pi^-$  is trivial on $U(p)$ and has highest weight  \[(a_0+p, 0,\dots , 0, -a_o-p ) \]   on $U(q)$ where $a_0= a-\frac{p+q-1}{2}$.  Hence the $K$-type  is trivial on the center of $K$ and has a 
$(K\cap K^\sigma ) $--fixed vector. 

The representation $\Pi^-$ is the  Langlands subrepresentation of a principal series representation $I(P,D, \nu)$ induced from a parabolic $P^- $= ${M^-} { A^-} { N^-}$, where  $A ^-$ is the vector subgroup of $U(p,q-2)$, of dimension $min(q-2,p-2)$. Its centralizer $M^-$ in $G$ is a smaller unitary group, which is non compact for $p,q\geq 3.$ (For details see \cite{KV} XI.10)

\medskip

{
\subsection{Relative Branching of $\Pi^- $ to $G'=U(p-1,q)$} \label{sec: relative branching }

We start by recalling the results of \cite{S}.



\medskip
\underline{Fact 1:} The restriction of the representation $\Pi^-$ to $G'\times U(1)$ is admissible, i.e. direct sum of irreducible representations, each occurring with multiplicity one.   
\cite{KO} \cite{S}.

\medskip
We refer to an irreducible  representation $\pi^-$ of $G'$ as {\em G'-type of  $\Pi^-$} if there is a character $\chi_r$ so that $\pi^- \otimes \chi_r$ is isomorphic to a summand  of $\Pi^-_{|G'} $. 

\medskip

We can use results about $K$-types to understand the restriction of  admissible representations due  to

\underline{Fact 2:} The restriction of the representation $\Pi^-$ to $G'\times U(1)$ is admissible implies the corresponding $(\bg,K)$-module is admissible. \cite{K1}

\medskip

Using pseudo dual pairs it is proved in  section 4 of \cite{S}  
that 

\underline{Fact 3:} The restriction of $\Pi_a^-$ to $G'_0 =G' \times U(1) $ is a direct sum of representations  $\pi_{\mu_{G'} }\otimes  \chi_r$. 
The $(\bg',K')$-module of the representation $\pi_{\mu_{G'}}$ of $G'$ is  cohomologically induced from  a $\theta$--stable parabolic $\bq'$   with Levi subgroup \[L' = L \cap G' =U(1) U(1) U(p,q-3) \]
and parameter
 \[ \mu_{G'}  =( a+1/2 +l,-a-1/2-k, 0, \dots 0) , \ \ \ \ \  l,k \in \bN.\]
 
 We denote the trivial character of $U(1)$ by $\chi_0$


 \begin{lemma} The  $G'$-types $\pi_{\mu_{G'}} $ of $\Pi_a^-$ are in the discrete spectrum of $G'/G' \cap H$ iff $\chi_0 \otimes \pi_{\mu_{G'}} $ has a nontrivial multiplicity in the restriction of $\Pi_a $ to $U(1) \times G' $. 
 \end{lemma}
 
 \proof Recall from \ref{sec:discrete spec} that $\pi_{\mu_{G'}}$ is in the discrete spectrum of  $L^2(G'/H') $ if k=l. 
 
  Consider the branching of $\Pi^-$ to $U(1)$. The subgroup $U(1)\subset U(p)$ commutes with $G'$ and hence acts by a scalar on the $(\bg',K')$-modules $A_{\bq'}(\mu_{G'}) $ . So it suffices to compute this scalar on the minimal $K$-type. In this case $K= U(p)U(q)$  and $K' = U(p-1)U(q)$. Using the formulas in  \cite{STV} we compute its action on $\A_{\bq'}(\mu_{G'}) $ by computing it on the factor $U(p)$ of its lowest $K'$-type. 
  
Recall that the minimal $K'$-type of a representation in the discrete spectrum of $Y^-$  is representation of $U(p-1) \times U(q)$  which is trivial on  the first factor.
  The highest weight of a representation of $U(p)$  whose restriction to $U(p-1)$ contains an $U(p-1)$-invariant vector
 is of the form $(r,0,0,... 0,-s)$ for positive integers $(r,s)$. By \cite{STV} $U(1)$ acts on its $U(p-1)$ invariant vector by the scalar $ l=r+s$. Hence 
 $A_{\bq'}(\mu_{G'}) $ is the $(\bg',K')$--module of a representation in the discrete spectrum of $Y^-$ iff $r=s$.
 \qed

  \medskip

 \subsubsection{Relative branching law}

   To complete the proof of the relative branching theorem for the restriction of $\Pi_a^-$ to $G'$  we have to show that  $\chi_0 \otimes \pi_{\mu_{G'}} $ has a nontrivial multiplicity in the restriction of $\Pi_a $ to $U(1) G' $ for $r=s \geq 0$. 
  We can now use $K$-type calculations or  branching by stages. We choose the latter since it complements well the considerations in section \ref{branchingPi+}. Recall the notation of \ref{sec:exhaustion} and \ref{sec:discrete spec}

 \begin{theorem} \label{theorem: Relative branchingPi-}
 
 Let  $\Pi_a^- $ be an irreducible representation in the discrete spectrum of G/H with nonsingular infinitesimal character and parameter $(0,\dots ,0,a, -a)$, $a\geq 0$.
\begin{itemize}
\item
  Let $\pi^-$ be an irreducible infinite dimensional representation with a nonsingular infinitesimal character in the discrete spectrum of $G'/H'$.  If
 \[ \mbox{Hom}_H(\Pi^-_{|G'} ,\pi^-)\not = 0\] then 
$\pi^-$ has a  parameter \[\mu_{G'}=( 0,\dots 0 ,a+k+1/2, -a-k-1/2) , \ \ \ \ \  k \in \bN\]
and \[ \mbox{dim Hom}_H(\Pi^-_{|G'} ,\pi^-) = 1\]

\item Conversely if $\pi^-$ is representation    with parameter
 \[ \mu_{G'}=(0,\dots 0, a+k+1/2, -a-k-1/2) , \ \ \ \ \  k \in \bN\]
Then $\chi_0\otimes  \pi^-  $  is a direct summand  of the representation $\Pi_a^-$ restricted to U(1)U(p-1,q).

\end{itemize}


\end{theorem}

\proof 
We already proved the first claim. It remains to show that that a  representations $\chi_0 \otimes \pi^- $ 
with parameter
 \[ \mu_{G'}=(0,\dots 0, a+k+1/2, -a-k-1/2) , \ \ \ \ \  k \in \bN\]
is a summand of $\Pi^-$ to $U(1)U(p-1,q)$.

As in the previous section we proceed by branching of  the representation $U^-(\ell)$ along the first sequence and second sequence
\begin{itemize}
\item $U(p-1,q) \subset O(2p-2, 2q) \times O(2) \subset O(2p, 2q)$
\item $U(p-1,q) \subset U(p,q) \subset O(2p,2q)$
\end{itemize}
Here we use the notation of \cite{K}. See also Section VI. and \ref{sec:exhaustion} of the definition of the parameter $\ell $.

We observe\\
\underline{Fact}  The representation $U^-(\ell)$ of $O(2p,2q)$ is isomorphic to $U^+(\ell)$ of $O(2q,2p)$. 

We  consider the restriction of the representation $U^+(\ell)$ of $O(2q,2p)$ to $O(2q,2p-2)\times O(0,2)$. By \ref{subsec:SO(2p,2q) to SO(0,2)SO(2p,2q-2)}
we conclude that it is 
a direct sum of representations $\chi_r^-  \otimes U^+(\ell + \text{\em k}) , k \in 2\bN$ where $\chi_r^-$ is a character of $O(2)$ which is nontrivial on the diagonal matrix (1,-1). (For the exact definition of $\ell$ see \cite{K}). So for $r=0$ and the using the isomorphism again we deduce that the restriction of $U^-(\ell )$ to $SO(2p-2,2q)$ is a direct sum of representations $U^-(\ell +k), k \in 2\bN$ .

Using the isomorphism again we see \ref{subsec:SO(2p,2q)-U(p,q)} that the discrete spectrum of the restriction of $U^+(\ell)$ of $O(2q,2p)$ to $U(q,p)$ contains exactly one representation $\Pi^+_{\ell/2,\ell/2}$, This representation is in  the discrete spectrum of $L^2(U(q,p)/U(1)U(q-1,p))$. Hence the representation $\Pi^+_{\ell/2,\ell/2}$ is the only representation in 
 the discrete spectrum of the restriction of $U^-(\ell)$ to $U(p,q-1)$
 which is $U(1)U(p,q-1)$ spherical.
 
 The same argument as in the previous section using the branching of $U^+(\ell)$ along the second sequence completes the proof.     
 \qed
}

\subsection{Restriction symmetry breaking } \label{sec: strong relative branching}

 \begin{prop}
Suppose that $b>a$ are positive integers and that the $(\bg,K)$ -module of 
$\Pi_a^-$ is a subset of $  C^\infty_{K}(X^-)\cap L^2(X^-)$.  Let 
 \[
 Res_{G'} :C^\infty_{K}(X^-) \rightarrow  C^\infty_{K'}(Y^-)\] be the restriction .
 Then the $(\bg',K')$--module of $\pi_b^-$ is a direct summand 
 of \[ Res_{G'} (\Pi_a ^- ).\]  
 
 \end{prop}
\proof
$\Pi^-_a$ is a direct sum of irreducible representation $\pi^-_{a+k} $ with $k\geq 0$. Recall that therefore its underlaying ($\bg,K$)-module  is a direct sum of the corresponding ($\bg',K'$)-modules . Recall also that the 
restriction 
\[Res_{G'} :C^\infty_{K}(X^-) \rightarrow  C^\infty_{K'}(Y^-)\] intertwines the $G'$ actions of the $K'$-finite functions.
So the image under restriction map of the $K$-finite functions in $\Pi^-_K$  is a sum of $(\bg',K')$- modules .

To prove the proposition it suffices to show that  $\pi^-_{a+k}$ is the $(\bg',K')$--modules of a summand in the restriction of $\Pi^-_a$ with parameter 
 $\mu_{G'}=(0,\dots 0, a+k+1/2, -a-k-1/2) , \  k \in \bN$ then $\pi^-_{a+k} $ is a submodule of 
$Res_{G'}(\Pi^-\cap C^\infty_{K}(X^-)) .$ To do this it suffices to show that if f is a K-finite function in $\Pi_a^-$ generating the $(\bg',K')$-submodule with parameter $\mu_{G'}$ then its restriction to $Y$ is non zero.

Let $\delta$ be the distribution on the functions in $C^\infty(X^-)$ at the identity coset
         \[\delta(f) = f(eU(1)U(p,q-1) )\]
and $\delta'$ the distribution on the functions in $C^\infty(Y^- )$ at the identity coset of $Y^-$
   \[\delta'(f) = f(eU(1)U(p-1,q-1) )\]   
  Since $Y^- \subset X^-$ we 
  get for $f \in \Pi^-\cap C^\infty_{K}(X^-)$
  \[ \delta(f) = \delta'(Res_{G'}(f)).\]
  Thus to prove the theorem it suffices to show that 
  if $f_{r,a_0}$ is a $K$-finite function in $\Pi^-_{a} \subset C^\infty_K(X^-)$ in a $K$-type with highest weight 
   \[(r,0,0,  ...,-r)\times (r+a_0,....,-a_0-r) \] which transforms under $U(p-1)$ by the trivial representation, i.e which is the lowest $K$-types of a $G'_0$- type of $\Pi^-$ 
      then 
  \[ \delta (f_{r,a_0})= f_{r,a_0}(eH') \not = 0 . \]
  
  Recall that $G= KBH^-$ and that $B$ is a one parameter group. Its centralizer $B^G$ is
  $U(p-1,q-1)$  and hence
  \[B^G \cap K = U(p-1) \times U(q-1). \]
  Furthermore
  \[B^G \cap G' = U(p-1,q-1)\]
  and 
  \[B^G \cap K' = U(p-1) \times U(q-1) .\]
  Thus $f_{r,a}$ is a $U(p-1)$ invariant function on $U(p)/U(p-1)$ and  an eigenfunction of the Laplacian and thus a  harmonic polynomial and  has nonzero constant term.

Alternatively, we may consider the operations of taking normal derivatives along $Y^-$ to implement the
symmetry-breaking operators.

 \qed

 \section{Closing remarks:} 
  
  As in the Gan-Gross-Prasad conjectures for discrete series representations the assumption, that the infinitesimal characters of the representations $\Pi^+, \Pi^- $ of G  are regular, is essential in our proofs and in the formulation of the main results in V.3.  Never the less a similar result maybe be true for representations also for representations with singular infinitesimal characters in the discrete spectrum of rank one symmetric spaces.
  
  \medskip
  
  We conjecture that the results  in V.3 depend only on the L-group of the symmetric space and hence are also true for representations in the discrete spectrum of symmetric spaces with the same L-group as  unitary symmetric spaces. This conjecture is supported by the results in \cite{KS} for orthogonal symmetric spaces.

 \bigskip
 
 \section{Appendix 1}
 Although we are making the assumption $ p,q\geq 3 $ we include a short discussion of the case p= 2, $q> 3$ to illustrate the difficulties in this case.
 
 An interesting example of relative branching is the restriction of a discrete series representation $\Pi_a $ of \[X=U(2,n)/U(1)U(1,n)\] to $G'=U(1,n)$.
{A discrete series representation of $X^+$ is also a discrete series representation of $U(2,n)$ since the Levi subgroup of the corresponding $\theta $--stable parabolic subgroup is compact \cite{KV}. 
 So we can directly apply the results of H.~He to determine the representations in  the discrete spectrum of its restriction to $U(2,n-1)$ \cite{H}. 
 Following H.~He we assign to the discrete series representation $\Pi_a^+$ the sequence 
 \[ (+,-, -  \dots, -, +)\] of 2+n signs showing that we have 2 noncompact simple roots for the positive roots defined by the Langlands parameter  of $\Pi^+_a$.
 A  discrete series representation $\pi $ of $U(1,n)$ has one signs $\oplus$ and 
 n-1 signs $\ominus$ . We have to align them according to the He rule:
 The only allowed combinations of 2 adjacent signs are 
 \[ (\oplus,+),(+\oplus) ,(-\ominus), (\ominus, -), (+,-) ,(-,+), (\oplus,\ominus),(\ominus,\oplus)\]
 
  Now suppose that we are restricting to the group $U(1,n)$. In this case  the only possible patterns for a discrete series representations of $U(1,n)$ which may appear in the restriction are
  \begin{itemize}
 \item one $\oplus $ in first place followed by n signs $\ominus$ 
 \item n signs $\ominus$ followed by a $\oplus$, 
  \end{itemize}
  These representations are either holomorphic or anti holomorphic representations.

 In the first case  we get one allowed pattern 
 \[+,\oplus ,\ominus,-,\ominus,-, \dots  ,\ominus , - ,+\]
 In the second case we have also one allowed pattern
 \[+,-,\ominus,-, \dots  ,\ominus ,-,\ominus,\oplus +\]
 To determine the infinitesimal characters of the discrete series representations in the discrete spectrum we have to considering the corresponding interlacing pattern for the infinitesimal characters.  
 We conclude that there are no discrete series representations $\pi$ in the discrete spectrum of the restriction of $\Pi_a^+$ to $U(1,n)$  with a He pattern in case 1 and 2 , since the corresponding interlacing pattern of the infinitesimal characters is not allowed. 
 Thus the restriction of $\Pi_a^+$ to $U(1,n)$ has no discrete spectrum \cite{GW}. 
  
   Observe that  
 the symmetric space $$Y^+=U(1,q)/U(1,0)U(q)$$ is a  Riemannian symmetric space  its discrete spectrum is empty, i.e. there are no non zero representations $\pi^+$. We conclude that the relative (strong) branching law is $(\Pi_a^+, \emptyset )$.

 Now consider the representations  $\Pi_a^- $ in the discrete spectrum of \linebreak $X^-=U(2,q)/U(1)U(2,q-1)$. They
  are no longer discrete series representations
 but they are cohomologically induced from a $ \theta $--stable parabolic with a noncompact 
 Levi subgroup $L^- = U(1)U(1)U(2,q-2)$.
 
  The symmetric space $$Y^-=U(1,q)/U(1)U(1,q-1) $$ is a rank one symmetric space of the group $U(1,q)$. It has a nonzero discrete spectrum  
 consisting of representations $\pi_b^-$ for b sufficiently large cohomologically induced from a $\theta $--stable parabolic subgroup  with a non compact Levi subgroup
 $ U(1)U(1)U(1,q-2)$.
 In addition here may be a finite number of singular representations. 
 The arguments in section VII show that the restriction of $\Pi_a$ to $U(1,q)$  is a direct sum of irreducible representations and  that representations $\pi_b^-$ are in the discrete spectrum of the restriction of $\Pi_a$ to $U(1,q)$ for $b>a$.


\bigskip

\section{Appendix 2} Though our main interest is in the case of unitary groups, we may
also treat in a similar way all rank one compact Riemannian symmetric spaces; their
spherical polynomials are again Jacobi polynomials $P_n^{(\alpha, \beta)}$. Thus in the quaternionic case
$$X_c = Sp(n+1)/Sp(1) \times Sp(n) = P^n(\IH)$$
with spherical polynomials
$$\varphi = P_n^{(2n-1, 1)}(x)$$
and also the octonionic case
$$X_c = F_4/Spin(9) = P^2(\mathbb O)$$
with spherical polynomials
$$\varphi = P_n^{(7,3)}(x)$$
and again $x = \cos(\theta)$ in terms of the polar angle (distance) $\theta$. For the
Flensted-Jensen discrete series we  consider in the quaternionic case the
subgroup $Sp(n)$ and the corresponding symmetric subspace \linebreak $X_c' = P^{n-1}(\IH)$.
The quaternionic hyperboloid is now
$$X = Sp(p,q+1)/Sp(p,q) \times Sp(1)$$
with the maximal compact subvariety $P^q(\IH)$. The Flensted-Jensen function
giving the lowest $K$-type is now again of product type
 $$\psi_{\lambda}(k b_s) = (\cosh s) ^{-i\lambda - \rho} P_n^{(2q-1,1)}(k)$$
for $k \in K,\, b_s \in B$. Here the parameters are (notation as in \cite{FJ}) 
\begin{eqnarray*}
i\lambda &= &2q-2p+1+n, (n \in 2\Z^+),\\  
\rho &= &2p+2q+1,\\ 
\rho_{\bt} &= &q+1, \\ 
\rho - 2\rho_{\bt} &= &2p-2q-1,\\
\mu_{\lambda} & = & i\lambda + \rho -2\rho_{\bt} = n.
\end{eqnarray*}
 The $n$ here is
even corresponding to spherical polynomials on \linebreak $X_c = P^q(\IH)$, and
the period integral is as before converging and non-zero for $n' \leq n$.
Note that since $i\lambda + \rho = 4q + n + 2$ we have the required decay,
since the radial density in this case equals (on $X_c = P^q(\IH)$)
$$D(s) =  (\cosh s)^{4q+3} (\sinh s)^{4p-1}.$$
And again there is a classical identity for these Jacobi polynomials
yielding the relative branching law via period integrals.


\end{document}